\documentclass[12pt]{article}
    \usepackage{srcltx}
   \usepackage{enumerate}
 \usepackage{latexsym}
   \usepackage{indentfirst}
\usepackage{amsmath,amsfonts,amssymb,amsthm}
 \usepackage[english]{babel}
\usepackage{graphicx}
\usepackage[font=small]{caption}
\usepackage{subfigure}
\usepackage[usenames]{color}
\theoremstyle{remark}
\newtheorem{rem}
{R~e~m~a~r~k~}

\newcommand{\new}{\newcommand}
\renewcommand{\textcolor}[1]{\relax}

\title{Universal kernel-type estimation of random fields
\thanks{
The study was supported by the program for fundamental
scientific research of the Siberian Branch of the Russian Academy of
Sciences, project FWNF-2022-0015. }}
\author{
Y.Y. Linke, I.S. Borisov,
and P.S. Ruzankin
\thanks{Sobolev Institute of Mathematics,
 Novosibirsk, 630090  Russia.}
\thanks{Novosibirsk State University, Novosibirsk, 630090  Russia.}
\thanks{E-mails: linke@math.nsc.ru, sibam@math.nsc.ru, ruzankin@math.nsc.ru.}
}

\date{}

\begin{document}
\new{\bb}{\mbox{\mathversion{bold}$\beta$}}
\new{\ee}{\mbox{\mathversion{bold}$\epsilon$}}
\new{\tha}{\mbox{\mathversion{bold}$\theta$}}
\new{\Tha}{\mbox{\mathversion{bold}$\Theta$}}
\new{\et}{\mbox{\mathversion{bold}$\eta$}}
\new{\ophi}{{\overline\varphi}}
\new{\oophi}{{\overline{\overline\varphi}}}
\new{\opsi}{{\overline\psi}}
\new{\pphi}{\mbox{\mathversion{bold}$\varphi$}}
\new{\ppsi}{\mbox{\mathversion{bold}$\psi$}}
\new{\rrho}{\mbox{\mathversion{bold}$\rho$}}
\new{\ddelta}{\mbox{\mathversion{bold}$\delta$}}
\new{\aalpha}{\mbox{\mathversion{bold}$\alpha$}}
\new{\xxi}{\mbox{\mathversion{bold}$\xi$}}

 \maketitle

 \abstract{Consistent weighted least square estimators are proposed for a wide class of nonpara\-metric regression models with random regression function, where
 this real-valued random function of $k$ arguments is assumed to be continuous with probability~1.
 We obtain explicit upper bounds for the rate of uniform convergence
 in probability of the new estimators to the unobservable random regression function
for both fixed  or random designs.
In contrast to the predecessors' results, the bounds
for the convergence are insensitive to the correlation structure of the $k$-variate design points.
As an application, we study the problem of estimating the mean and covariance functions of random fields with additive noise  under dense data conditions.
\textcolor{blue}{The theoretical results of the study are illustrated by simulation examples which show that the new estimators are more accurate in some cases than the Nadaraya--Watson ones.
 An example of processing real data
on earthquakes in Japan in 2012--2021 is included. }

}

 \begin{center}
\vspace{0.5cm}

{\bf 1. Introduction}
\end{center}
We study the following
regression model:
 \begin{equation}  \label{f2}
 Y_i=f(X_i)+\xi_i,
 \qquad i=1,\ldots,n,
 \end{equation}
where $f(t)$,  $t:=(t^{(1)},\ldots,t^{(k)})\in \Theta\subset {\mathbb R}_+^k$, $k\ge 1$, is an unknown real-valued random function (a random field). We assume that $\Theta$ is a compact set,
the random field $f(t)$ is continuous  on $\Theta$ with probability 1, and
the design $\{X_i;\,i=1,\ldots,n\}$ consists of a collection of  observed random vectors with unknown   (generally speaking) distributions,  not necessarily independent or identi\-cally distributed.
The random design points $X_i$ may depend on $n$, i.e.,
a triangular array scheme for the design can be considered within this model.
In particular, this scheme includes regression models with fixed design.
Moreover, we do not require that the random field $f(t)$ be
independent of the design $\{X_i\}$.

Next, we will assume that the unobservable random errors $\{\xi_i\}$ (a noise) form a martingale difference sequence, with
\begin{equation}  \label{xi}
M_p:=\sup_i{\bf E}|\xi_i|^p<\infty \quad\mbox{for some} \,\, p>k \,\,\mbox{and}\,\, p\ge 2.
\end{equation}
We also assume that $\{\xi_i\}$ are independent of the collection $\{X_i\}$ and the random field $f(\cdot)$. The noise $\{\xi_i\}$ may depend on $n$.

Our goal is to construct consistent in $C(\Theta)$ estimators for the random regression field $f(t)$ by the observations $\{(X_i,Y_i);\, i\le n\}$
under minimal restrictions on the design points $\{X_i\}$, where $C(\Theta)$ denotes the space of all continuous functions on $\Theta$ with the uniform norm.

 In the classical case of nonrandom  $f(\cdot)$, the most popular estimating procedures are based on kernel-type estimators. We emphasize among them the Nadaraya--Watson,
the Priestley-Chao,  the Gasser-M\"uller estimators, the local polynomial estimators, and some of their modifications  (see
H\"ardle, 1990;  Wand and Jones,  1995; Fan and  Gijbels,  1996;
\textcolor{blue}{ Fan and  Yao, 2003;
  Loader, 1999;
 Young, 2017};    M\"{u}ller, 1988).
\textcolor{blue}{
We do not aspire for
providing a comprehensive review of this actively
developing (especially in the last two decades)
area of nonparametric estimation, and will focus only on publications representing certain methodological areas.
We are primarily interested in conditions on the design elements. In this regard, a large number of publications in this area may be tentatively
divided into the two groups:
the studies dealing with fixed design $\{X_i;\, i\le n\}$  or with random one.
In papers dealing with a random design, as a rule,
the design consists of  independent identically distributed random variables
  or stationary observations satisfying known forms of dependence, e.g.,
various types of mixing conditions,
association, Markov or martingale properties, etc. Not attempting to present a comprehensive review, we may note the papers by
}
Kulik and Lorek (2011), Kulik and Wichelhaus (2011), Roussas (1990, 1991), Gy\"orfi et. al. (2002), Masry (2005), Hansen (2008), Honda (2010), Laib and Louani (2010), \textcolor{blue}{Li et al. (2016), Hong and  Linton (2016), Shen and  Xie  (2013),    Jiang and Mack (2002),    Linton and   Jacho-Chavez  (2010), Chu and  Deng (2003) (see also the references in the papers).
Besides, in the recent studies by Gao et al. (2015), Wang and   Chan (2014),  Chan and  Wang (2014), Linton and Wang  (2016), Wang and Phillips (2009a,b),  Karlsen et al. (2007), the authors considered nonstationary design sequences under special forms of dependence (Markov chains, autoregressions, sums of moving averages, and so on).
}

\textcolor{blue}{In the case of fixed design, the vast majority of papers make certain assumptions on the regularity of design (see
  Zhou and Zhu, 2020; Benelmadani et al., 2020; Tang et al., 2018; Gu et al., 2007;  Benhenni et al., 2010;
M{\"u}ller and Prewitt, 1993; Ahmad and Lin, 1984; Georgiev, 1988, 1990). In univariate models, the nonrandom design points $X_i$ are most often restricted by the formula $X_i=g(i/n)+o(1/n)$ with a function $g$ of bounded variation, where the error term $o(1/n)$ is uniform in $i=1,\ldots,n$. If $g$ is linear, then the design is equidistant. Another regularity condition in the univariate case is the relation $\max\nolimits_{i\leq n}(X_i-X_{i-1})=O(1/n)$, where the design elements are arranged in increasing order.
In a number of recent studies, a more general condition $\max\nolimits_{ i\leq n}(X_i-X_{i-1})\to 0$
can be found
(e.g., see Yang and Yang, 2016; He, 2019; Wu et al., 2020). In several works, including
those dealing with the so-called weighted estimators, certain conditions are imposed on the behavior of functions of design elements, but meaningful corresponding examples
are limited to cases of regular design (e.g., see Zhang et al., 2019 ; Zhang et al., 2018; Liang and Jing, 2005; Roussas et al., 1992; Georgiev, 1988).}

     The problem of uniform approximation of the kernel-type 
estimators has been studied by many authors
(e.g., see Einmahl and Mason, 2005; Hansen, 2008; \textcolor{blue}{ Gu et al., 2007;  Shen and   Xie, 2013; Li et al., 2016; Liang and  Jing,  2005; Wang and  Chan, 2014;  Chan and   Wang, 2014; Gao et al., 2015}  and the references therein).

In this paper, we study a class of kernel-type estimators, asymptotic properties of
which do not depend on the design correlation structure. The design may be fixed (and not necessarily regularly spaced) or random
(with not necessarily weakly dependent components).
We present weighted least square estimators where the weights are chosen as the Lebesgue measures of the elements of a finite random partition
of the  regression function domain $\Theta$
such that every partition element corresponds to one design point. As a result, the proposed kernel estimators for the regression function are
transformation of sums of weighted observations in a certain way with the structure of multiple Riemann integral sums, so that conceptually our approach is close to
the methods of Priestley and Chao (1972) and of Mack and M{\"u}ller (1988), who considered the cases of univariate fixed design and i.i.d. random design, respectively.
Explicit upper bounds are obtained for the rate of uniform convergence of these estimators to the random regression field.

In contrast to the predecessors' results, we do not impose any restrictions on the design correlation structure.
\textcolor{blue}{
We will consider
the maximum cell diameter of the above-mentioned
partition of $\Theta$
generated by the design elements,
 as the main characteristic of the design.
Sufficient conditions for the consistency of the new estimators, as well as the windows' widths
will be derived in terms of that characteristic.
The advantage of that characteristic over the classical weak dependence conditions is that
the characteristic
 is insensitive to forms of correlation of the design elements. The main condition will be that the maximum cell diameter tends to zero in probability as the sample size grows.
Note that such requirement is, in fact, necessary, since only when the design densely fills the regression function domain, it is possible to reconstruct the function more or less precisely. }

\textcolor{blue}{Univariate versions of this estimation problem were studied in Borisov et al. (2021) and Linke et al. (2022) where the asymptotic analysis and
simulations showed that the proposed estimators perform better than
the Nadaraya--Watson ones in several cases.
Note that the univariate case in Borisov et al. (2021)
does not allow direct generalization to a multivariate case,
since the weights were defined  there as the spacings of the variational series generated by the design elements.
Note also that the estimator
in Borisov et al. (2021)
is a particular univariate case of the estimators proposed in this paper, but not the only one.
One of the univariate estimators studied here  
may be more accurate than the estimator in Borisov et al. (2021) (see Remark~\ref{r3} below).
Conditions on the design elements similar to those of this paper were used
Linke and Borisov (2022), and in Linke (2023). The conditions provide uniform consistency of the estimators,
but guarantee only pointwise consistency of the Nadarya--Watson ones.
Besides, similar restrictions on the design elements were used before in Linke and Borisov~(2017, 2018), and Linke (2019) in estimation of the parameters of several nonlinear regression models.}

\textcolor{blue}{In this paper,
we will assume that the unknown random regression function $f(t)$, $t\in \Theta$, is continuous with probability 1.
Considering the general case of random regression function allows us to obtain results on estimating the mean function of a random regression process. In regard to estimating random regression functions, we may note the papers by
Li and Hsing (2010), Hall et al. (2006), Zhou et al. (2018), Zhang and  Wang (2016, 2018),  Yao  et al. (2005),  Zhang and  Chen (2007),   Yao (2007),
Lin and  Wang (2022). In those papers, the mean and covariance functions of the random regression process~$f(t)$ were estimated when, for independent noisy copies of the random process, each of the trajectories
was observed in a certain subset of design elements (nonuniform random time grid).
Estimation of mean and covariance of random processes is an actively developing area of nonparametric estimation, especially in the last couple of decades, is of independent interest, and plays an important role in subsequent analyses
(e.g., see Hsing and Eubank, 2015; Li and Hsing, 2010;  Zhang and Wang, 2016; Wang et al., 2016).
 }

\textcolor{blue}{
Estimation of random regression functions usually deals with
either random
or deterministic design.
In the case of random design, it is usually assumed that the design elements are independent identically distributed (e.g., see Hall et al., 2006; Li and Hsing, 2010; Zhou et al., 2018; Yao , 2007; Yao et al., 2005; Zhang and Chen, 2007; Zhang and Wang, 2016, 2018; Lin and Wang, 2022). Some authors emphasized that their results can be extended to weakly dependent design (e.g., see Hall et al., 2006). For deterministic time grids, regularity conditions are often required (e.g., see Song et al., 2014; Hall et al., 2006).
In regard to denseness of filling the regression function domain, the two types of design are distinguished in the literature: either the design is ``sparse'', e.g., the number of design elements in each series is uniformly limited (Hall et al., 2006; Zhou et al., 2018; Li and Hsing, 2010), or the design is ``dense'' and the number of elements in a series increases with the sequential number of the series (Zhou et al., 2018; Li and Hsing, 2010).
Uniform consistency of several estimators of the mean of random regression function
was considered, for example, by Yao et al. (2005), Zhou et al. (2018), Li and Hsing (2010), Hsing and Eubank (2015), Zhang and Wang (2016).
}

\textcolor{blue}{
In this paper, we consider one of the variants of estimation of the mean of a random regression function as an application of the main result.
 In the case of dense design,
 uniformly consistent estimators are constructed for the mean function, when the series-to-series-independent design is arbitrarily correlated
inside each series. We require only that, in each series, the design elements form a refining partition of the domain of
the random regression function.
Our settings also include a general deterministic design situation, but we do not impose traditionally used regularity conditions.
Thus, in the problem of estimating the mean function, as well as in the problem of estimating the function in
model~(\ref{f2}),
 we weaken  traditional conditions on the 
 design elements.
Note that methodologies used for estimating the mean function for dense and for sparse data usually differ
(e.g., see Wang et al., 2016). In the case of growing number of observations in each series, it is natural to preliminarily
evaluate the trajectory of the random regression function in each series and then average the estimates over all series (e.g., see Hall et al., 2006).
That is what we will do in this paper following this generally accepted approach.
}
Universal estimates both for the mean  and covariance functions of a random process in the case of sparse data, insensitive to the nature of the dependence of design elements, are proposed in Linke and Borisov (2024).

This paper has the following structure. Section~2 contains the main results on the rate of uniform convergence of the new estimators to the random regression function. In Section~3, we consider an application of the main results to the problem of estimating the mean and covariance function of a random regression field. In Section~4, the asymptotic normality of the new estimators is discussed.
\textcolor{blue}{Section~5 contains several simulation examples. In Section~6, we discuss an example of assessing real data on earthquakes in Japan in 2012--2021.
In Section 7, we summarize the main results of the paper.
The proofs of the theorems and lemmas from Sections 2--4 are contained in Section 8.
}

\newpage

\begin{center}
{\bf 2. Main assumptions and results}
\\
\end{center}

Without loss of generality we will assume that $d(\Theta\cup 0)\le 1$, where $$d(A):=\sup_{x,y\in A}\|x-y\|$$
is the diameter of a set $A$ and $\|\cdot\|$ is the
supnorm in ${\mathbb R}^k$.
In what follows, unless otherwise stated, all the limits will be taken as $n\to\infty$.

Our approach recalls a construction
  of multivariate Riemann integrals. To this end, we need  the following condition on the design $\{X_i\}$.

$({\bf D})$ {\it For each $n$, there exists a random partition of the set $\Theta$ into $n$
Borel-measurable subsets $\{\Delta_i; \;i=1,\ldots,n\}$ such that
$\delta_n:=\max_{i\le n}d(\Delta_i\cup X_i)\to 0$ in probability.
}

\textcolor{black}{
Condition (D) means that, for every  $n$, the set $\{X_i;\,i\le n\}$ forms a $\delta_n$-net in the compact set $\Theta$. In particular, Condition (D) is satisfied if the design points $\{X_i\}$ are pairwise distinct, $X_i\in \Delta_i$ for all $i\le n$, and
$\max_{i\le n}d(\Delta_i)\to 0$ in probability.
}

In the case $\Theta=[0,1]^k$,
a regularly spaced design satisfies Condition (D).
Moreover, if $\{X_i;\,i\ge 1\}$ is a stationary sequence satisfying an $\alpha$-mixing condition and $[0,1]^k$ is the support
 of the distribution of $X_1$, then Condition $(D)$ is fulfilled (see Remark~3 in Linke and Borisov, 2017).
 It is not hard to verify that, for i.i.d. design points with the probability density function of $X_1$  bounded away from zero on $[0,1]^k$, one can have $\delta_n=O\left(\frac{\log n}{n^{1/k}}\right)$ with probability 1. Notice that the dependence of random variables $\{X_i\}$  in Condition $(D)$ may be much stronger than that in these examples
(see Linke and Borisov, 2017, 2018 \textcolor{blue}{and the example below}).

\textcolor{blue}{
{\it Example}. Let a sequence of bivariate random variables
$\{{X}_i;\,i\ge 1\}$ is defined by the relation
\begin{equation}\label{w13}
{ X}_{i}=\nu_{i}{ U}_{1i}+(1-\nu_{i}){ U}_{2i},
\end{equation}
 where the random vectors $\{{ U}_{1i}\}$ and  $\{{U}_{2i}\}$ are independent
 and uniformly distributed on the rectangles $[0,1/2]\times [0,1]$ and $[1/2,1]\times[0,1]$, respectively,
 while the sequence
 $\{\nu_i\}$ does not depend on  $\{{ U}_{1i}\}$ and  $\{{U}_{2i}\}$ and consists of Bernoulli random variables with success probability $1/2$,
 i.e., the distribution of  $X_i$ is the equi-weighted mixture of the two above-mentioned uniform distributions.
The dependence between the random variables $\{\nu_i\}$  is defined by the equalities
$\nu_{2i-1}=\nu_1$ and $\nu_{2i}=1-\nu_1$. In this case, the random variables $\{{ X}_i;\,i\ge 1\}$ in (\ref{w13}) form a stationary sequence of random variables uniformly distributed on the unit square $[0,1]\times[0 ,1]$,
 but, say, all known mixing conditions are not satisfied here because,
for all natural $m$ and $n$,
\begin{eqnarray*}
{\mathbb P}\big({ X}_{2m}\in [0,1/2]\times [0,1],\,{X}_{2n-1}\in [0,1/2]\times [0,1]\big)=0.
\end{eqnarray*}
On the other hand, it is easy to check that the stationary sequence $\{{X}_i\}$ satisfies the Glivenko--Cantelli theorem. This means that, for any fixed $h>0$,
$$\#\{i:\,\|{ t}-{ X}_{i}\|\le h,\, 1\le i\le n\}\sim 4h^2 n$$ almost surely uniformly in  $t$, where $\#$ denotes the standard counting measure. In other words, the sequence
$\{{ X}_i\}$ satisfies Condition  $(D)$.
}

\textcolor{blue}{
It is clear that, according to the scheme of this example, one can construct various sequences of dependent random variables uniformly distributed on $[0,1]\times [0,1]$,
based on the choice of different sequences of the Bernoulli switches with the conditions $\nu_{j_k}=1$ and $\nu_{l_k}=0$ for infinitely many
indices $\{j_k\}$ and $\{l_k\}$, respectively. In this case, Condition $(D)$ will also be satisfied.}
\textcolor{blue}{
 But the corresponding sequence $\{{ X}_i\}$ (not necessarily stationary) may not satisfy the strong law of large numbers. For example, a similar situation occurs when $\nu_j=1-\nu_1$ for $j=2^{2k-1},\ldots,2^{2k}-1$ and $\nu_j=\nu_1$ for $j=2^{2k},\ldots,2^{2k+1}-1$, where $k=1,2,\ldots$
(i.e., we randomly choose one of the two rectangles $ [0,1/2]\times [0,1] $ and $ [1/2,1] \times [0,1] $,
into which we randomly throw the first point, and then alternate the selection of one of the two rectangles
by the following numbers of elements of the sequence: $1$, $ 2$, $2^2$, $2^3$, etc.).
Indeed, we can introduce the notation
$n_k=2^{2k}-1$,  $\tilde n_k=2^{2k+1}-1$, $S_m=\sum\nolimits_{i=1}^m X_{i}^{(1)}$, with ${ X}_i=(X_{i}^{(1)}, X_{i}^{(2)})$, and note
 that, for all outcomes consisting the event  $\{\nu_1=1\}$, one has
\begin{equation}\label{example2}
\frac{S_{n_k}}{n_k}=\frac{1}{n_k}\sum\limits_{i\in N_{1,k}}U_{1i}^{(1)}+\frac{1}{n_k}\sum\limits_{i\in N_{2,k}}U_{2i}^{(1)},
\end{equation}
}\textcolor{blue}{
where ${ U}_{ji}=(U_{ji}^{(1)}, U_{ji}^{(2)})$, $j=1,2$;   $N_{1,k}$ and $N_{2,k}$  are the collections of indices for which the observations
$\{X_i, i\leq n_k\}$ lie in the rectangles $[0,1/2]\times [0,1]$ or $[1/2,1]\times [0,1]$, respectively.
 It is easy to see that $\#(N_{1,k})=n_k/3$ and   $\#(N_{2,k})=2\#(N_{1,k})$. Hence,   ${S_{n_k}}/{n_k}\to{7}/{12}$ almost surely as  $k\to\infty$
 due to the strong law of large numbers for the sequences  $\{U_{1i}^{(1)}\}$ and  $\{u_{2i}^{(1)}\}$. On the other hand,
 for all elementary outcomes in the event  $\{\nu_1=1\}$,   as $k\to\infty$, we have with probability 1
 \begin{equation}\label{example2+}
\frac{S_{\tilde n_k}}{\tilde n_k}=\frac{1}{\tilde n_k}\sum\limits_{i\in \tilde N_{1,k}}U_{1i}^{(1)}+\frac{1}{\tilde n_k}\sum\limits_{i\in \tilde N_{2,k}}U_{2i}^{(1)}\to \frac{5}{12},
\end{equation}
where  $\tilde N_{1,k}$ and $\tilde N_{2,k}$  are the collections of indices for which the observations
$\{{ X}_i, i\leq \tilde n_k\}$ lie the rectangles $[0,1/2]\times [0,1]$ or  $[1/2,1]\times [0,1]$, respectively.
In proving the convergence
in (\ref{example2+}) we took into account that $\#(\tilde N_{1,k})=(2^{2k+2}-1)/3$, $\#(\tilde N_{2,k})=2n_k/3$, i.e.,  $\#(\tilde N_{1,k})=2\#(\tilde N_{2,k})+1$.
}

\textcolor{blue}{
Similar arguments are valid for elementary outcomes consisting the event  $\{\nu_1=0\}$.
$\hfill\Box$
}

In what follows,  by $K(s)$,  $s\in \mathbb R^k$, we will denote the kernel function. We assume that
the kernel function is zero outside $[-1,1]^k$ and is a centrally symmetric probability density function, i.e.,
$K(s)\ge 0 $, $K(s)=K(-s)$ for all $s\in [-1,1]^k$, and $\int_{[-1,1]^k} K(s) ds =1$.
For example, we may consider product-kernels of the form
$$K(s)=\prod\limits_{j=1}^kK_o(s^{(j)}),
$$
where $K_o(\cdot)$  is a univariate symmetric probability density function with support $[-1,1]$.
We also assume that the function $K(s)$ satisfies the Lipschitz condition with  constant $L\ge 1$:
\begin{equation*}\label{lip}
|K(x)-K(y)|\le L\left(|x^{(1)}-y^{(1)}|+\cdots+|x^{(k)}-y^{(k)}|\right)
\end{equation*}
for all $x=(x^{(1)},\dots,x^{(k)})$ and $y=(y^{(1)},\dots,y^{(k)})$,
and put $K(y)=0$ for all $y$ such that $\|y\|> 1$.
Notice that, under these restrictions, $\sup_s K(s)\le L$.

Put
$$K_{\varepsilon}(s):=\varepsilon^{-k} K(\varepsilon^{-1}s).$$
By $\theta_{\varepsilon}$ we denote a random vector  with the density $K_{\varepsilon}(t)$,
which is independent of the random variables $\{Y_i\}$.

Let $\Lambda(\cdot)$ denote the Lebesgue measure in ${\mathbb R}^k$.
Introduce the following notation:
\begin{eqnarray}\label{main}
f^*_{n,\varepsilon}(t):=\frac{\sum_{i=1}^nY_{i}K_{\varepsilon}(t-X_{i})\Lambda(\Delta_i)}{\sum_{i=1}^nK_{\varepsilon}(t-X_{i})\Lambda(\Delta_i)},
\end{eqnarray}
where $0/0=0$  by definition;
\begin{eqnarray}\label{smooth}
J_{\varepsilon}(t):=\int\limits_{\Theta}K_{\varepsilon}(t-x)\,\Lambda(dx)\equiv
{\bf P}(t-\theta_{\varepsilon}\in {\Theta}),\quad t\in \Theta;
\end{eqnarray}
$$\omega_f(\varepsilon):=\sup_{x,y\in \Theta:\, \|x-y\|\le
\varepsilon}|f(x)-f(y)|.
$$

Now, notice that

$${f^*_{n,\varepsilon}}(t)={\rm arg}\min\limits_{z\in \mathbb R}
\sum\limits^n_{i=1}(Y_{i}-z)^2K_{\varepsilon}(t-X_{i})\Lambda(\Delta_i),
$$
i.e., the estimators of the form
(\ref{main})
belong to the class of {\it weighted least square estimators}. 
Estimators (\ref{main})
are also called {\it local constant } ones.

Finally,  we will assume that there exist constants $\rho>0$ and $\varepsilon_0\in (0,1]$ such that
\begin{equation}\label{rho}
J_{\varepsilon}(t)\ge \rho\ \mbox{ for all }\ t\in\Theta\,\  \mbox{and positive}\,\,\varepsilon \le\varepsilon_0.
\end{equation}
So, some cases (for example, if
$\Theta$ contains isolated points) are excluded from the scheme under consideration.

 \begin{rem}
 Notice that if the set $\Theta$ can be represented as the union of hyperrectangles with the edges of lengths greater than $\varepsilon_0$ and kernel $K$ is a product-kernel, then we have the lower bound $\rho \ge 2^{-k}$.
\end{rem}

The main result of this paper is as follows.

{\bf Theorem 1}. {\it Let the conditions $(\rm D)$ and  $(\ref{rho})$ hold. Then, for any fixed $\varepsilon\in (0,\varepsilon_0]$,
\begin{eqnarray}\label{T12}
\sup_{t\in \Theta}|f^*_{n,\varepsilon}(t)-f(t)|\le \omega_f(\varepsilon)
+\zeta_n(\varepsilon)
\end{eqnarray}
 with probability $1$, where $\zeta_n(\varepsilon)$ is a positive random variable such that
\begin{eqnarray}
{\bf P}\left( \zeta_n(\varepsilon)>y\right)&\le&
G(k,p)\, \rho^{-p}\, M_p\, L^{p/2}\, y^{-p}
\,\varepsilon^{-k(p/2+1)}\,{\bf E}(\delta_n^{kp/2})\,
\notag \\
&&
\label{rest2}
+\ {\bf P}(\delta_n>\varepsilon
\min\{1,\, \rho(k2^{k+1} L)^{-1}\}),
\end{eqnarray}
where
$$0<G(k,p)< (p-1)^{p/2} 2^{p(k+(3/2))}
\,
\left( 1
+\frac{k}{2^{(p-k)/(p+1)} - 1}
\right)^{p+1}.
$$
}

In what follows, we will denote by
    $O_p(\eta_n)$ some  univariate random variables $\zeta_n$ such that, for all $y>0$,
   $$\limsup_{n\to\infty} {\bf P}(|\zeta_n|/\eta_n>y)\le \beta(y),
   $$
where  $\{\eta_n\}$ is a sequence of positive nonrandom numbers, $\lim_{y\to\infty}\beta(y)=0$, and the function $\beta(y)$ does not depend on $n$.

\begin{rem}
 For example, let the function $f$ be nonrandom.  In
 (\ref{rest2}), put  $$y=\left(\varepsilon^{-k(p/2+1)}\,{\bf E}(\delta_n^{kp/2})\right)^{1/p}.$$ Applying the power Markov's inequality with exponent  ${kp/2}$ for the second summand in  (\ref{rest2}), we obtain that, under the conditions of the theorem,
$$ \zeta_n(\varepsilon)=O_p\left(\left(\varepsilon^{-k(p/2+1)}\,{\bf E}(\delta_n^{kp/2})\right)^{1/p}\right)
$$
and  there exists a solution $\varepsilon\equiv\varepsilon(n)$ to the equation
\begin{eqnarray}\label{opt}
{\bf E}(\delta_n^{kp/2})=\varepsilon^{k(p/2+1)}\omega^p_f(\varepsilon).
\end{eqnarray}
It is clear that this solution vanishes as $n\to \infty$.
In fact, the value $\varepsilon(n)$ minimizes in $\varepsilon$ the order of smallness for the right-hand side of (\ref{T12}).
Notice that   $\delta_n/\varepsilon(n)\stackrel{p}{\to} 0$ and $(\varepsilon(n))^{-k(p/2+1)}\,{\bf E}(\delta_n^{kp/2})$ in view of (\ref{opt}).
\end{rem}

\textcolor{black}{
Taking Remark 1 into account one can obtain the following two assertions as consequences of Theorem $1$.
}

{\bf Corollary 1}. {\it
Let $\cal C$ be a set of nonrandom equicontinuous functions from the function space
$C[0,1]^k$ $($for example, a precompact subset of $C[0,1]^k$$)$. Then,
under Condition $(D)$,
$$
\gamma_n({\cal C}):=\sup_{f\in \cal C}\,\sup_{t\in [0,1]^k}| f^*_{n,\tilde\varepsilon(n)}(t)-f(t)|\stackrel{p}\to 0,
$$
where $\tilde\varepsilon(n)$ is a solution to  equation $(\ref{opt})$ in which the modulus of continuity
$\omega_{f}(\varepsilon)$ is replaced with the universal modulus
$\omega_{\cal C}(\varepsilon):=\sup_{f\in \cal C}\omega_f(\varepsilon)$.
In this case, $\gamma_n({\cal C})=O_p(\omega_{\cal C}(\tilde\varepsilon(n)))$}.

{\bf Corollary 2}. {\it If the modulus of continuity of the regression random field $f(t)$, $t\in[0,1]^k$, in  Model $(\ref {f2})$ meets the condition
$\omega_{f}(\varepsilon)\le \zeta \tau(\varepsilon)$ a.s., where $\zeta>0$ is a proper random variable and $\tau (\varepsilon)$ is a positive continuous nonrandom function, with
 $\tau (\varepsilon)\to 0$ as $\varepsilon\to 0$, then,
 under Condition $(D)$,
\begin{eqnarray}\label{optim3}
\sup_{t\in [0,1]^k}| f^*_{n,\hat\varepsilon(n)}(t)-f(t)|\stackrel{p}\to 0,
\end{eqnarray}
where $\hat\varepsilon(n)$ is a solution to equation $(\ref{opt})$ in which the modulus of continuity $\omega_f(\varepsilon)$ is replaced
with  $\tau (\varepsilon)$.}

{\it Example} 2. Let $\Theta=[0,1]^k$, $\delta_n\le \nu n^{-1/k}$, with ${\bf E}\nu^{kp/2}<\infty$, and $\omega_f(\varepsilon)\le \zeta \varepsilon^{\alpha}$, $\alpha\in (0,1]$, where
$\zeta$ is a proper random variable. Then $\varepsilon(n)=O\left(n^{-\frac{1}{2k(1/p+1/2)+\alpha}}\right)$
and
$$\sup_{t\in [0,1]^k}|f^*_{n,\varepsilon}(t)-f(t)|=O_p\left(n^{-\frac{\alpha}{2k(1/p+1/2)+\alpha}}
\right).
$$

In particular, in the one-dimensional case, if $f(\cdot)=W(\cdot)$ is a Wiener process on $[0,1]$, and the i.i.d. random variables $\xi_i$ are centered Gaussian, then
by L\'{e}vy's modulus of continuity theorem,
for any arbitrarily small $\nu>0$, we have
$$\sup_{t\in [0,1]}|f^*_{n,\varepsilon}(t)-f(t)|=O_p(n^{-1/3+\nu}).
$$
Here we put  $k=1$, $\alpha=1/2-\nu_1$, and $1/p<\nu_2$, with arbitrarily small positive $\nu_1$ and $\nu_2$.

\begin{rem}\label{r3}
\textcolor{blue}{
Let $k=1$, $\Theta=[0,1]$.  Denote by
$X_{n:1}\leq \ldots\leq X_{n:n}$ the ordered  sample 
$\{X_i;\,i=1,\ldots, n\}$.
Put
 $$X_{n:0}:=0,\quad X_{n:n+1}:=1,\quad  \Delta_{ni}:=(X_{n:i-1},\, X_{n:i}], \quad i=1,\ldots,n.$$
Denote by $Y_{ni}$
the response variable $Y$  corresponding to $X_{n:i}$ in  (\ref{f2}).
Then we can write down estimator (\ref{main}) for the function $f$ as
\begin{equation}\label{est4}
 f^*_{n,\varepsilon}(t)=\frac{\sum_{i=1}^nY_{ni}K_{\varepsilon}(t-X_{n:i})
\Delta X_{ni}}{\sum_{i=1}^nK_{\varepsilon}(t-X_{n:i})\Delta X_{ni}},
\end{equation}
where
$$ \Delta X_{ni}:= \Lambda(\Delta_{ni}) = X_{n:i}-X_{n:i-1}.$$
This estimator was proposed and  studied in detail in Borisov et al. (2021).}

\textcolor{blue}{
But, 
instead of $\{\Delta_{ni}\}$, we can consider Voronoi cells
$$ \widetilde\Delta_{ni}:= \left( \frac{X_{n:i-1}+X_{n:i}}{2}, \,
\frac{X_{n:i}+X_{n:i+1}}{2} \right]$$}
\textcolor{blue}{
and write down the corresponding estimator:
}
\textcolor{blue}{
\begin{equation}\label{est4+}
 \widetilde f^*_{n,\varepsilon}(t)=\frac{\sum_{i=1}^nY_{ni}K_{\varepsilon}(t-X_{n:i})
\widetilde\Delta X_{ni}}{\sum_{i=1}^nK_{\varepsilon}(t-X_{n:i})
\widetilde\Delta X_{ni}},
\end{equation}
where
$$ \widetilde\Delta X_{ni}:= \Lambda(\widetilde\Delta_{ni}) =
\frac{X_{n:i+1}-X_{n:i-1}}{2}. $$
Repeating, for the last estimator, the corresponding proofs in Borisov et al. (2021) originally
applied to estimator (\ref{est4}), we can easily see that all properties
of estimator (\ref{est4}) are retained for (\ref{est4+}), except
the constant factor in the asymptotic variance.
Namely, let the regression function $ f (t) $ 
be twice continuously differentiable and nonrandom, let the errors $ \{\xi_i\} $ be
independent identically distributed, centered with finite second moment $M_2$, and independent of the design
 $\{X_i\}$,
whose elements be independent identically distributed. In addition, let $ X_1 $ have
 a strictly positive density $ p(t) $ which is continuously differentiable. 
Then
\begin{equation}\label{var}
{\mathbb Var} f_{n,\varepsilon}^*(t)\sim \frac{2M_2}{h np(t)}\int\limits_{-1}^1K^2(u)du,\qquad
{\mathbb Var}\widetilde f_{n,\varepsilon}^*(t)\sim \frac{1.5M_2}{h np(t)}\int\limits_{-1}^1K^2(u)du.
\end{equation}
The former asymptotic relation was established in Lemma 3 by Borisov et al. (2021).
The latter relation can be proved by repeating the proof of that lemma with obvious changes.
Hence, in the case of independent and identically distributed design points, the asymptotic variance of the estimator can be reduced by choosing an appropriate partition.
}

\textcolor{blue}{
Thus, for $k=1$, this paper deals with a more general class of estimators (\ref{main}) than
that in Borisov et al. (2021) where estimator (\ref{est4}) was studied, and representatives of the new class can have certain advantages over the estimator (\ref{est4}).}
\end{rem}

\newpage

\begin{center}
{\bf 3. Application to estimating the mean and covariance functions \\ of a random regression function}
\end{center}

In this section, as an application of Theorem 1, we will construct a consistent estimator for the mean function of the random  regression function in Model $(\ref {f2})$.
We consider the following multivariate statement of the  problem of estimating the mean function of an a.s. continuous random regression stochastic process. Consider $N$ independent copies of Model $(\ref {f2})$:
 \begin{equation}  \label{f3}
 Y_{i,j}=f_j(X_{i,j})+\xi_{i,j},
 \qquad i=1,\ldots,n,\,\,\,\,j=1,\ldots,N,
 \end{equation}
where $f(t), f_1(t),\ldots, f_N(t)$, $t\in [0,1]^k$, are i.i.d. unknown a.s. continuous stochastic processes, and, for every $j$, the collection $\{\xi_{i,j};\,i\le n\}$ satisfies  condition (\ref{xi}).
Here and in what follows, the subscript $j$
denotes the sequential
number of such a copy.
Introduce the notation
$$\overline{f^*}_{N,n,\varepsilon}(t):=\frac{1}{N}\sum\limits_{j=1}^Nf^*_{n,\varepsilon,j}(t).
$$


{\bf Theorem 2}. {\it Let Condition $(D)$ for Model $(\ref{f2})$ be fulfilled and
\begin{equation}\label{supmoment}
{\bf E}\sup_{t\in [0,1]^k}|f(t)|<\infty.
\end{equation}
Besides, let a sequences $\varepsilon\equiv \varepsilon_n\to 0$ and a sequence of naturals $N\equiv N_n\to \infty$ satisfy the conditions
\begin{equation}\label{Neps}
\varepsilon^{-k(p/2+1)}\,{\bf E}(\delta_n^{kp/2})\to 0\,\,\,\mbox{and}\,\,\,N{\bf P}(\delta_n>\varepsilon
\min\{1,\, \rho(k2^{k+1} L)^{-1}\})\to 0.
\end{equation}
Then
\begin{equation}\label{LLN}
\sup\limits_{t\in [0,1]^k}|\overline{f^*}_{N,n,\varepsilon}(t)-{\bf E}f(t)|\stackrel{p}\to 0.
\end{equation}
}

\begin{rem}
If condition (\ref{supmoment}) is replaced with a slightly stronger condition $${\bf E}\sup_{t\in [0,1]^k}f^2(t)<\infty$$ then, under the restrictions (\ref{Neps}), one can prove the uniform consistency of the estimator
$$M_{N,n}^*(t_1,t_2):=\frac{1}{N}\sum\limits_{j=1}^Nf^*_{n,\varepsilon,j}(t_1)f^*_{n,\varepsilon,j}(t_2),\,\,\,\,t_1,t_2\in [0,1]^k,
$$
for the unknown mixed second-moment function ${\bf E}f(t_1)f(t_2)$, where $\varepsilon\equiv \varepsilon_n$ and $N\equiv N_n$ are defined in (\ref{Neps}). The proof is based on the same arguments as those in proving Theorem 2, and therefore is omitted.
In other words, under the above-mentioned conditions, the estimator $${\rm Cov}^*_n(t_1,t_2):= M_{N,n}^*(t_1,t_2)-\overline{f^*}_{N,n,\varepsilon}(t_1)\overline{f^*}_{N,n,\varepsilon}(t_2)$$ is uniformly consistent for the covariance function of the random regression field $f(t)$.
\end{rem}

\newpage
\begin{center}
{\bf 4. Asymptotic normality}
\end{center}
 In this section, we discuss sufficient conditions for  asymptotic normality of the estimators $f^*_{n,\varepsilon}(t)$.
Denote by ${\cal F}_0$ the trivial $\sigma$-field, and by ${\cal F}_j$  the $\sigma$-field generated by the collection $\{\xi_1,\ldots,\xi_j\}$, by the design points, and by the regression random field. 

{\bf Theorem 3}. {\it
Let the design $\{X_i\}$ do not  depend on $n$.  
Under Condition $(D)$, assume that, for some $t\in \Theta$ and a sequence $\varepsilon\equiv \varepsilon_n$,
\begin{equation}\label{AN1}
h_n:=\frac{\max_{j\le n} \big(K_{\varepsilon}\left(t-X_j\right)
\Lambda(\Delta_j)\big)^2}{\sum_{j=1}^n \big(K_{\varepsilon}\left(t-X_j\right)
\Lambda(\Delta_j)\big)^2}\stackrel{p}\to 0,
\end{equation}
$${\bf E}(\xi_j^2\  |\ {\cal F}_{j-1})=\sigma^2\ \mbox{ a.s. for all }j,
$$
$$ \max_j {\bf E}\left(\xi_j^2\, {\bf 1} (\xi_j^2>a/h_n)\  |\  {\cal F}_{j-1}\right)
\stackrel{p}{\to} 0\ \mbox{ for all }a>0.
$$
Then
\begin{equation*}
B^{-1}_{n,\varepsilon}(t)\left(f^*_{n,\varepsilon}(t)-f(t) -r_{n,\varepsilon}(t)\right)\stackrel{d}\to N(0,\sigma^2),
\end{equation*}
where $N(0,\sigma^2)$ is a centered Gaussian random variable with variance $\sigma^2$,
$$B^2_{n,\varepsilon}(t):=J^{-2}_{n,\varepsilon}(t)
\sum\limits_{i=1}^n\big(K_{\varepsilon}(t-X_{i})\Lambda(\Delta_i)\big)^2,
$$
$$r_{n,\varepsilon}(t):=J^{-1}_{n,\varepsilon}(t)
\sum\limits_{i=1}^n(f(X_{i})-f(t))K_{\varepsilon}\left(t-X_{i}\right)\Lambda(\Delta_i),
$$
$$J_{n,\varepsilon}(t):=\sum\limits_{i=1}^nK_{\varepsilon}(t-X_{i})\Lambda(\Delta_i).
$$

}

The theorem is a direct consequence of Corollary 3.1 in Hall and Heyde (1980).

\begin{center}
{\bf 5. Simulation examples}
\end{center}

In this section, we present simulations comparing the estimator $f^*_{n,\varepsilon}(t)$
defined in (\ref{main})
with the Nadaraya--Watson estimator
$$
\hat f_{n,\varepsilon}(t):=
\frac{\sum_{i=1}^nY_{i}K_{\varepsilon}(t-X_{i})}{\sum_{i=1}^nK_{\varepsilon}(t-X_{i})}
$$
in the 2-dimensional case.
For this estimator, we will assume $0/0=0$, like that was done for the estimator  $(\ref{main})$.

The elements of the design space $\Theta = [-1,1]\times [-1,1]$ will be denoted by $(x,y)$.
The following two algorithms were used to partition the space $\Theta$ into the sets $\Delta_i$.

The first algorithm is the Voronoi partitioning.
For each $i$, the set $\Delta_i$ is the Voronoi cell corresponding to $X_i$, i.e. the set of all points of $\Theta$ that lie closer to $X_i$ than to any other design point. The {\sf deldir} R package was employed for calculation of the squares of the cells.

The second algorithm is recursive partitioning  by coordinate-wise medians.
First, we divide $\Theta$ into the two rectangles by the line
$t^{(1)} = {\rm median} \{X_1^{(1)},\dots,X_n^{(1)} \},$
where the median is the midpoint of the interval
$\left(X_{\lfloor n/2 \rfloor}^{(1)},\  X_{\lfloor n/2 \rfloor+1}^{(1)} \right)$
when all the points are sorted
in increasing order with respect to the first coordinate.
Then each of the two rectangles is divided recursively. If, at some step,
a rectangle contains two or more design points then it is divided into the two parts:
If the rectangle's width is greater than its height, then the rectangle
is divided by
the line $t^{(1)}= {\rm median} \{X_{l}^{(1)}: l\in  B\},$ where $B$ is the set of indices
of the design points falling into the rectangle;
otherwise it is divided  by
the line $t^{(2)}= {\rm median} \{X_{l}^{(2)}: l\in  B\}$.
As soon as there is only one design point $X_i$ in a rectangle, the rectangle is
put to be $\Delta_i$.

\begin{figure}[!h]
    \centering
    \subfigure[]{
        \includegraphics[width=0.45\textwidth]{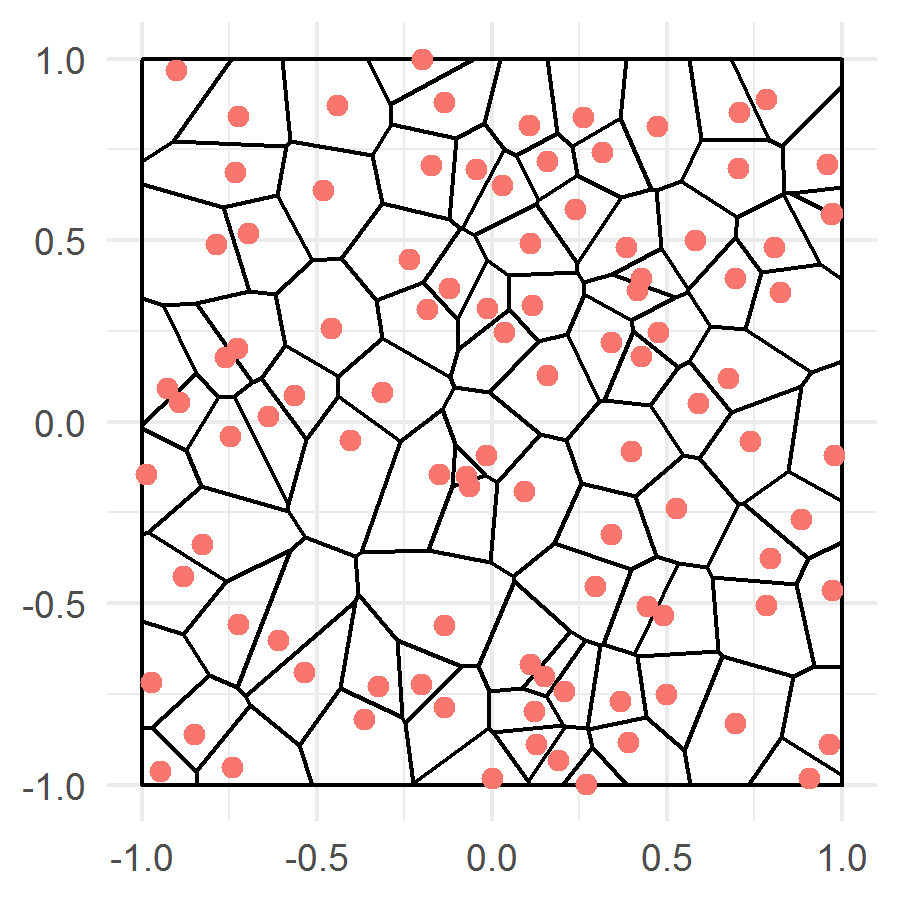}
    }
    \subfigure[]{
        \includegraphics[width=0.45\textwidth]{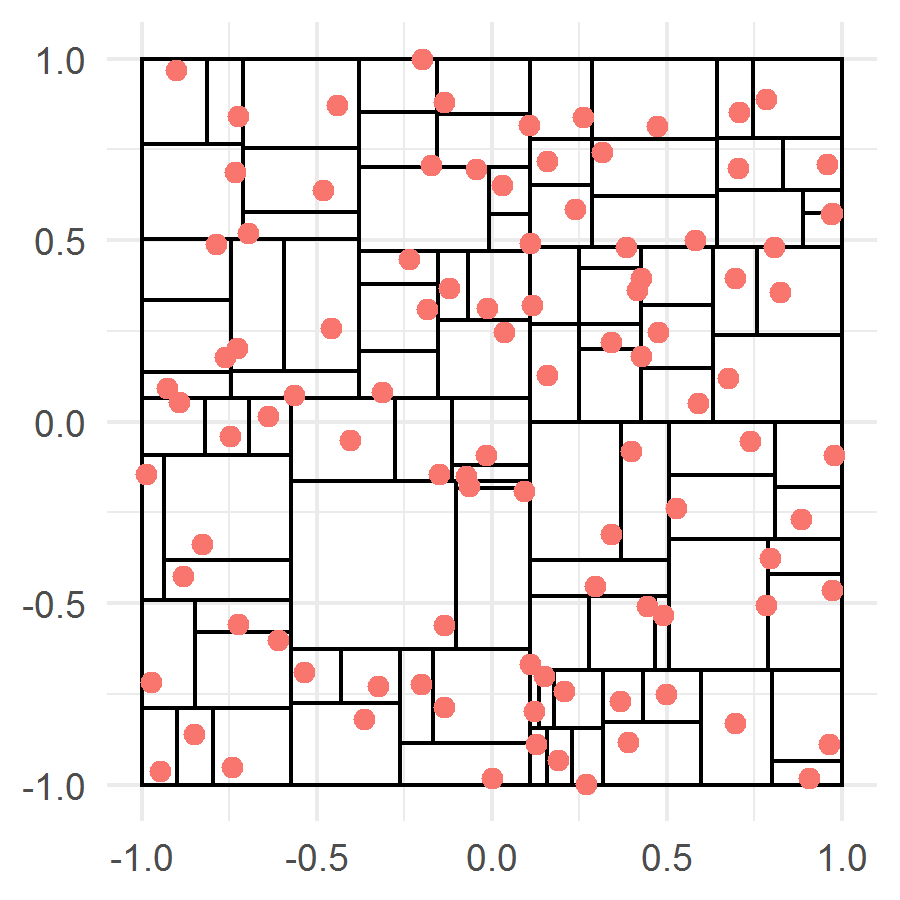}
    }
    \vspace{-12pt}
    \caption{Partitioning into Voronoi cells (left) and partitioning by coordinate-wise medians (right) for the same collection of points.}  \label{partfig}
\end{figure}

Results of partitioning $\Theta$ into cells for a collection of 100 points by the both algorithms are displayed in Fig.~\ref{partfig}.

In the simulation examples below, we used the tricubic kernel
$$K(x,y)=\frac{440}{162\pi} \max\left\{0, \left(1-\sqrt{x^2+y^2}^{\,3}\right)^3\right\}.$$
In each example, 1000 simulation runs were performed. In each of the simulation runs, 5000 design points were generated and randomly divided into the training (80\%) and validation (20\%) sets.
For the design points $X_i$, the observations $Y_i$ were generated with i.i.d. Gaussian noise with standard deviation $\sigma=0.5$.
For each of the tested algorithms, on the training set, the optimal $\varepsilon$ was calculated by 10-fold cross-validation minimizing the average of mean-square errors.
The $\varepsilon$ was selected from 20 values located on the logarithmic grid from 0.01 to 0.5.
The random partitioning for the cross-validation was the same for all the tested algorithms.

Then, for each of the algorithms, the model, trained on the training set with the chosen $\varepsilon$, was used to compute the mean-square error (MSE) for the observations of the validation set:
$$ {\rm MSE} = \frac{1}{m} \sum_j (f^*_{\varepsilon}(X_j)-Y_j)^2,$$
where the sum is taken over the validation set, and $m$ is the size of the set.
Besides, that model was employed to compute the maximal absolute error (MaxE) for the true values of the target function $f$ on the $100\times 100$ uniform lattice on $\Theta$:
$$ {\rm MaxE} = \max_j |f^*_{\varepsilon}(\gamma_j)-f(\gamma_j)|,$$
where the maximum is computed for the elements $\gamma_j$ of the $100\times 100$ lattice covering $\Theta$.

The algorithms that were compared will be denoted by NW (Nadaraya-Watson), ULCV (Universal Local Constant estimator (\ref{main}) with Voronoi partitioning), and ULCM (Universal Local Constant estimator (\ref{main}) with coordinate-wise Medians partitioning).

The results of the simulation runs are presented as median (1-st quartile, 3-rd quartile) and are compared
between the estimators with the paired Wilcoxon test.

In the examples below, we intentionally chose the densities of the design points with high nonuniformity in order to demonstrate possible advantages of the new estimator.

\bigskip
\noindent\textbf{5.1. Example 1}

\begin{figure}[!h]
    \centering
    \subfigure[Design points] {
        \includegraphics[width=0.35\textwidth]{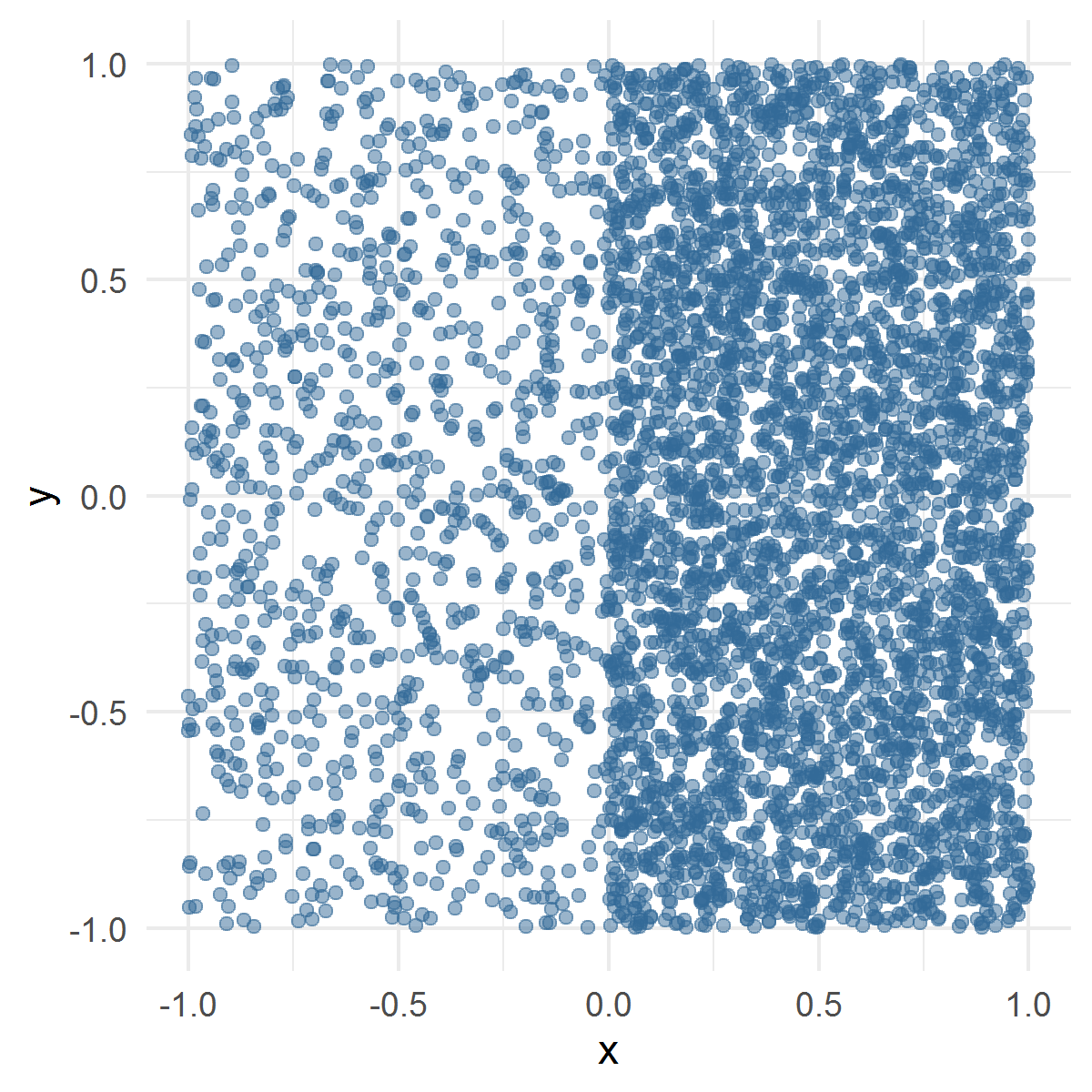}
    }
    \subfigure[MSE]{
        \includegraphics[width=0.25\textwidth]{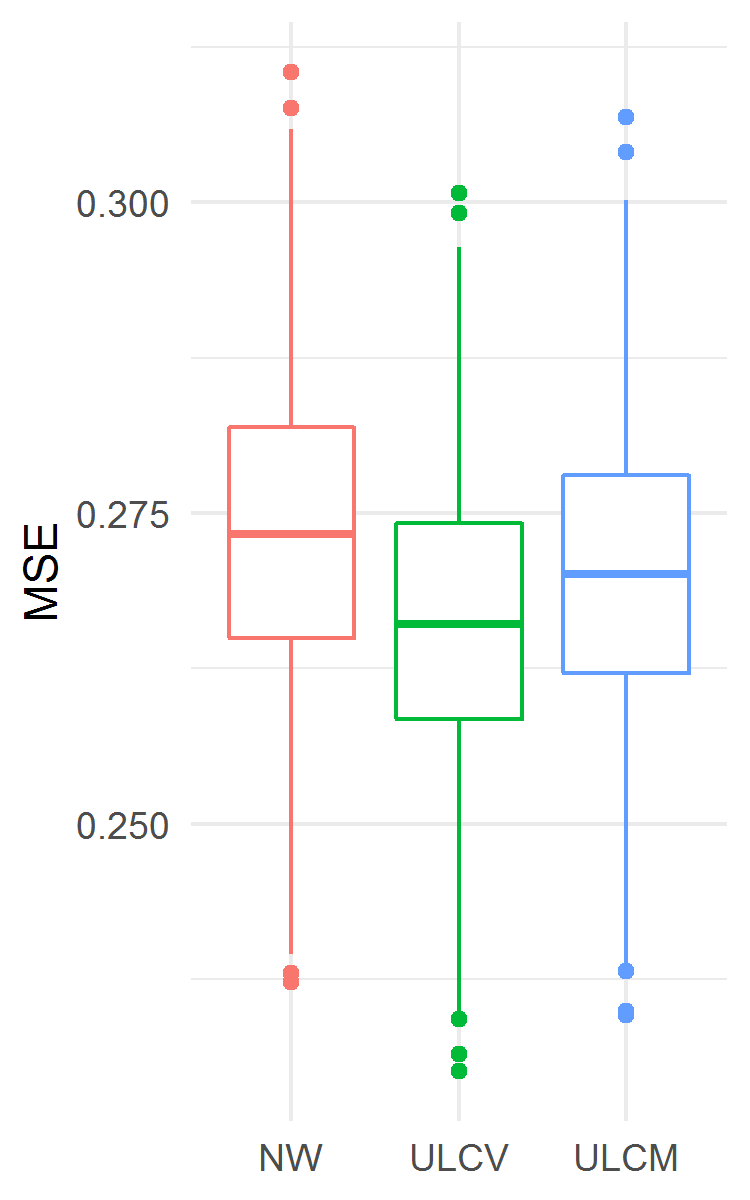}
    }
    \subfigure[MaxE]{
        \includegraphics[width=0.25\textwidth]{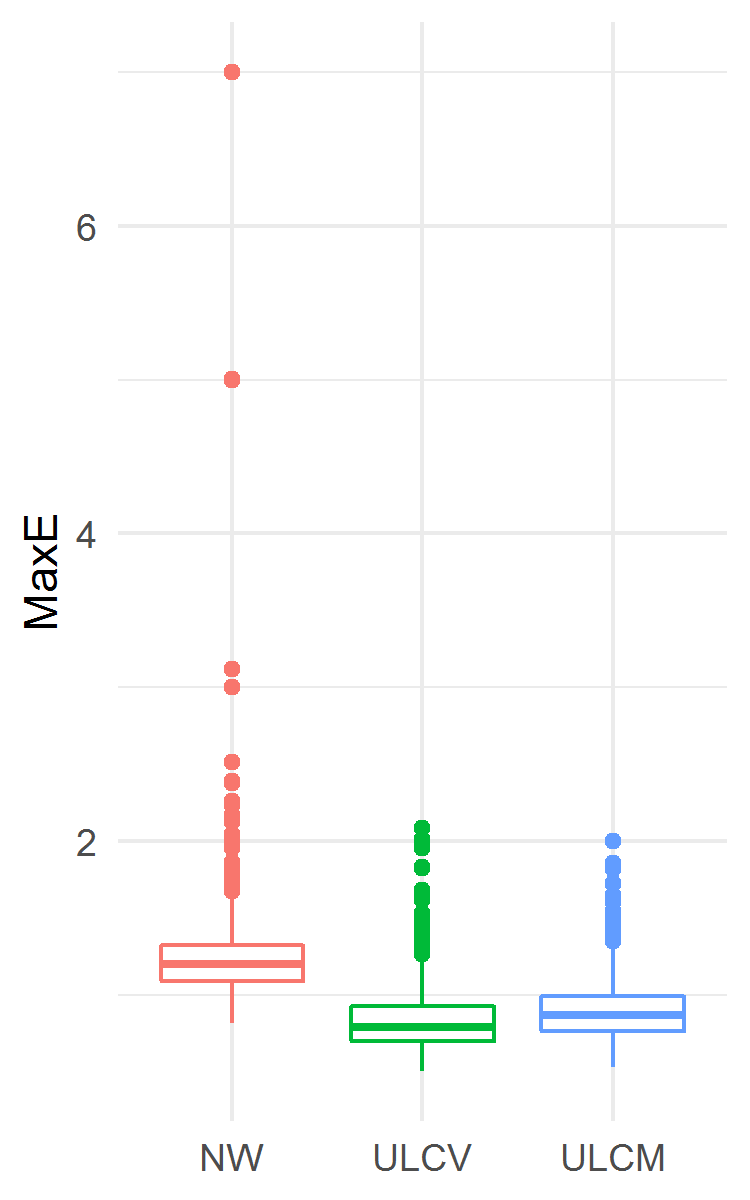}
    }
    \vspace{-12pt}
    \caption{Design points in an Example~1 experiment (left), the mean-square errors (middle), and maximal absolute errors (right) in Example~1}\label{ex1mse}
\end{figure}

\begin{figure}[!h]
    \centering
    \subfigure[The estimated function $z=f(x,y)$]{
        \includegraphics[width=0.45\textwidth]{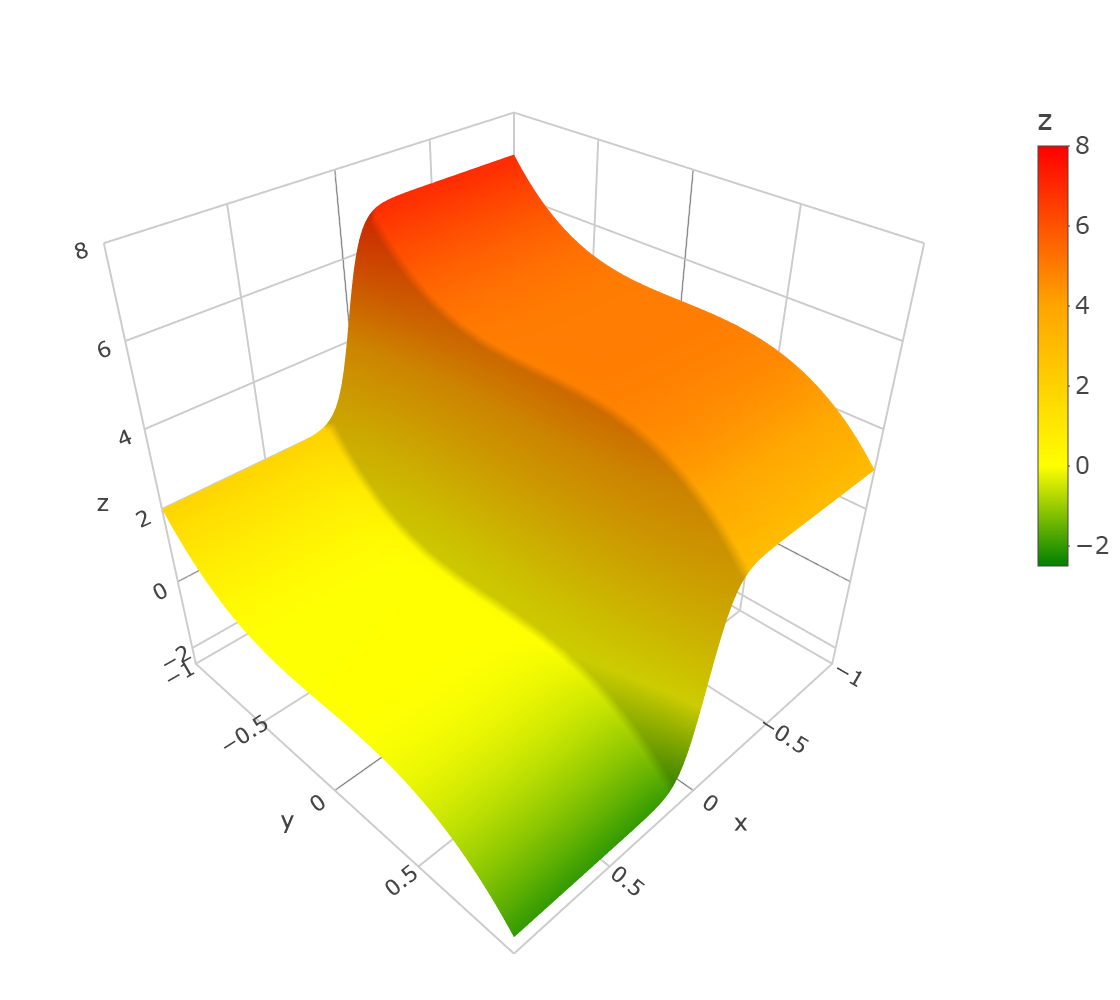}
    }
    \subfigure[A ULCV estimate $z=f^*(x,y)$]{
        \includegraphics[width=0.45\textwidth]{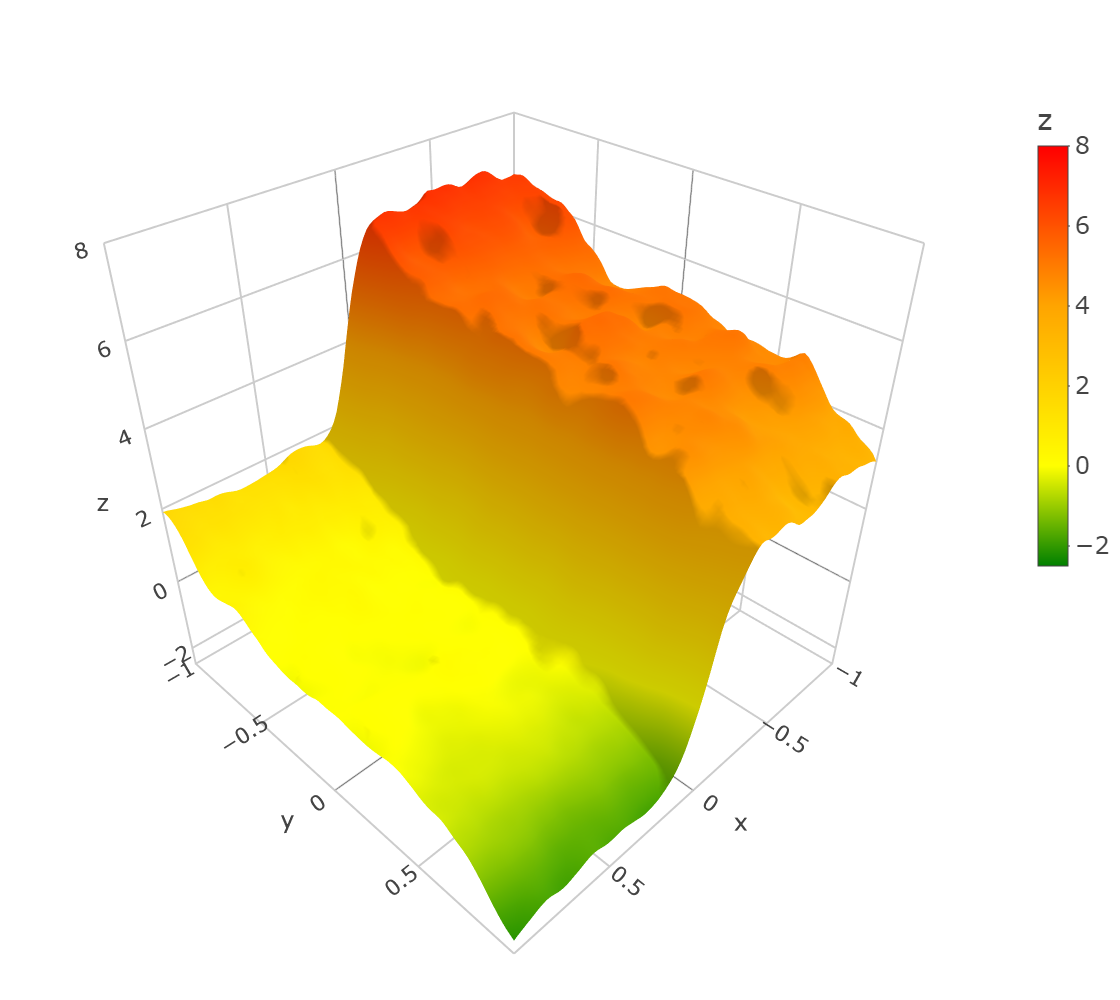}
    }
    \vspace{-12pt}
    \caption{The estimated function (left) and a result of ULCV estimator (right) in Example~1}\label{ex1fun}
\end{figure}

In this example, we approximate the nonrandom regression function
$$f(x,y)=\frac{5}{1+e^{-20x}} - 2y^3.$$
The  design points were generated in a way similar to that in Example in Sec.~1.
First, we choose the left rectangle $[-1,0)\times [-1,1]$ or the right rectangle $[0,1]\times [-1,1]$
with equal probabilities and draw $X_1$ uniformly distributed in the chosen rectangle. Then we draw $10$ design points
uniformly distributed in the other rectangle. Then we draw $10^2$ design points uniformly distributed in the rectangle where $X_1$ lies.
Then we draw $10^3$ design points uniformly distributed in the rectangle where $X_2$ lies, and so on.
In other words, we alternate the rectangle after $1$, $10$, $10^2$, $10^3$, ... draws.
One draw of the $5000$ design points is depicted in Fig.~\ref{ex1mse}. The estimated function $f(x,y)$ and a computed ULCV estimate are depicted in Fig.~\ref{ex1fun}.

%
%

The results are presented in Fig.~\ref{ex1mse}. The ULCV estimator appeared to perform best among the three considered ones both for MSE and MaxE accuracy measures. In particular, the ULCV estimator was better than the NW one: MSE 0.2661 (0.2584,  0.2742) vs. 0.2734 (0.2650, 0.2819), $p<0.0001$; MaxE 0.7878 (0.7013,  0.9230) vs. 1.1998 (1.0911, 1.3250), $p<0.0001$.
In this example, the ULCM estimator was better than the NW one as well.

\bigskip


\noindent\textbf{5.2. Example 2}

In this example, we approximate the nonrandom regression function
\begin{equation}\label{funex2}
f(x,y)=\sin\left(10 \sqrt{x^2+y^2}\right)/\sqrt{x^2+y^2}.
\end{equation}

\begin{figure}[!h]
    \centering
    \subfigure[Design points] {
        \includegraphics[width=0.35\textwidth]{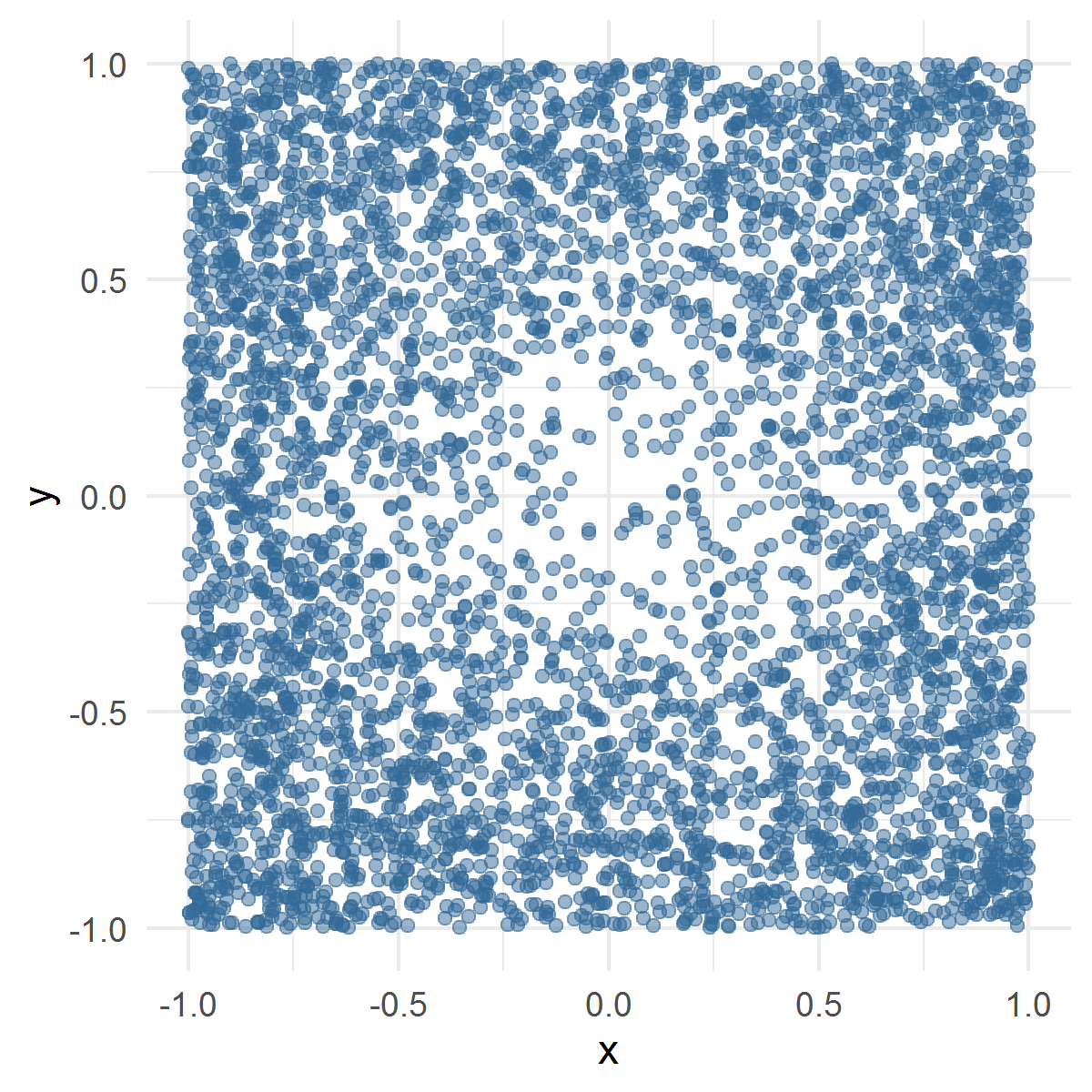}
    }
    \subfigure[MSE]{
        \includegraphics[width=0.25\textwidth]{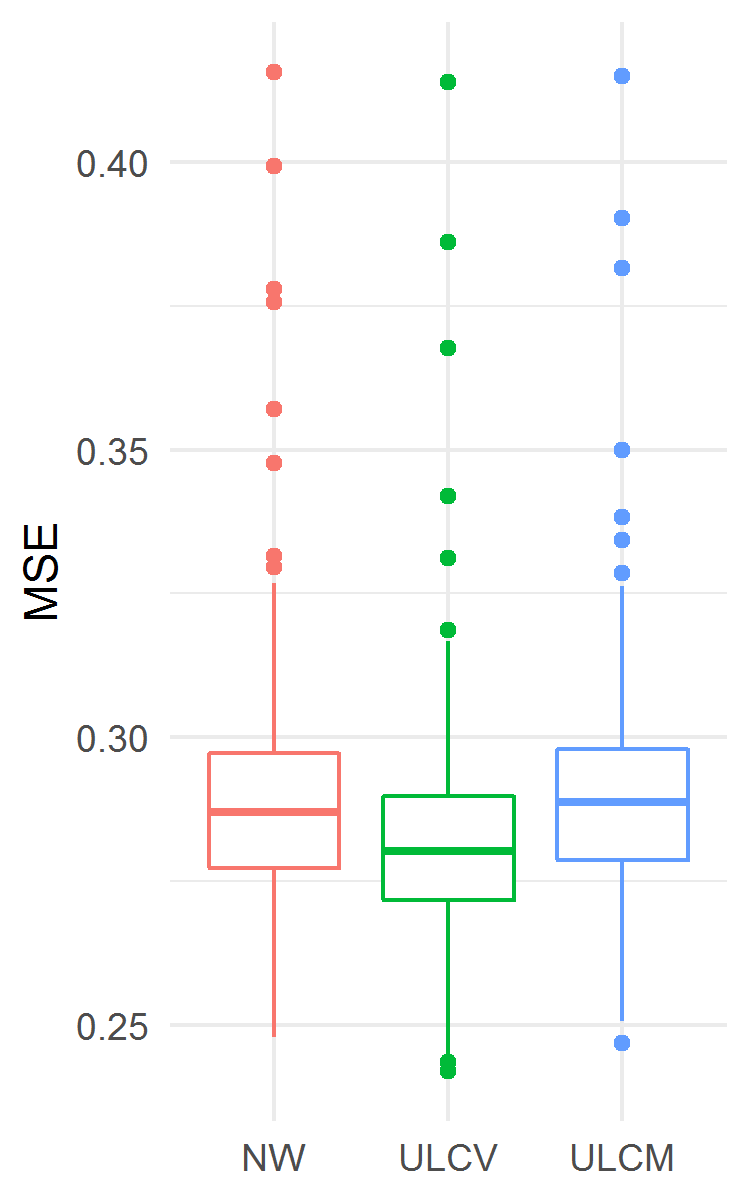}
    }
    \subfigure[MaxE]{
        \includegraphics[width=0.25\textwidth]{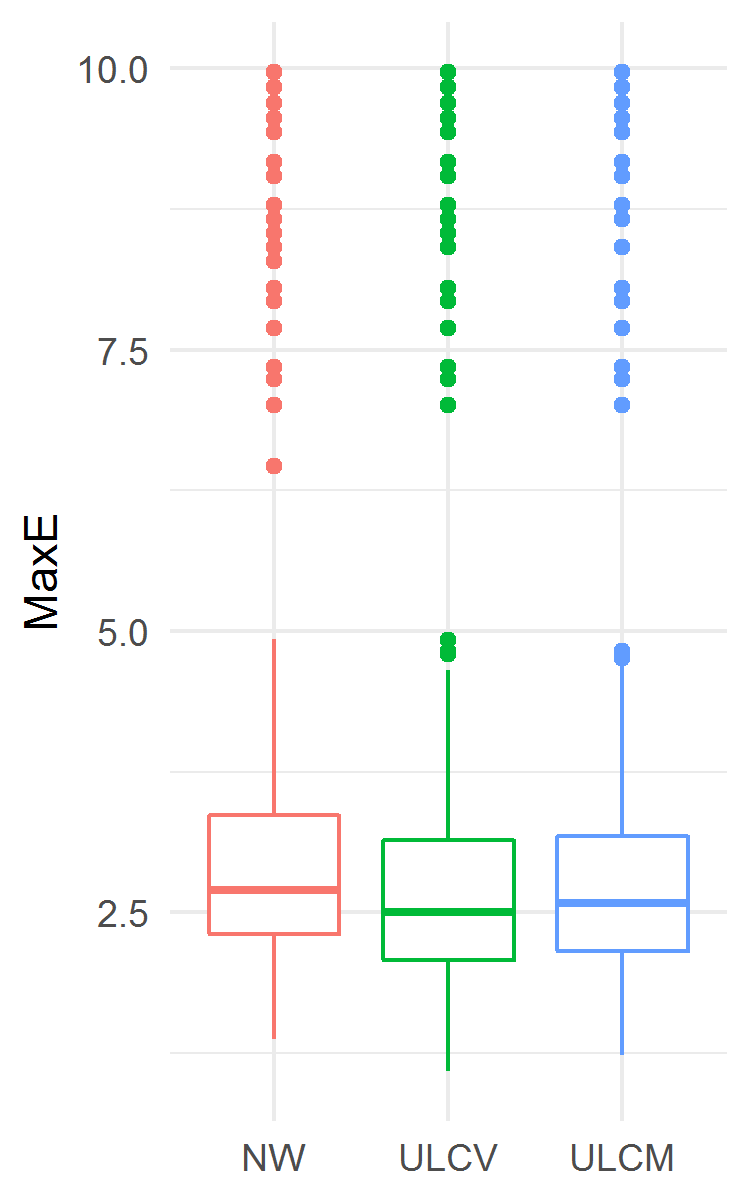}
    }
    \vspace{-12pt}
    \caption{Design points in an Example~2 experiment (left), the mean-square errors (middle), and maximal absolute errors (right) in Example~2}\label{ex2mse}
\end{figure}

\begin{figure}[!h]
    \centering
    \subfigure[The estimated function $z=f(x,y)$]{
        \includegraphics[width=0.45\textwidth]{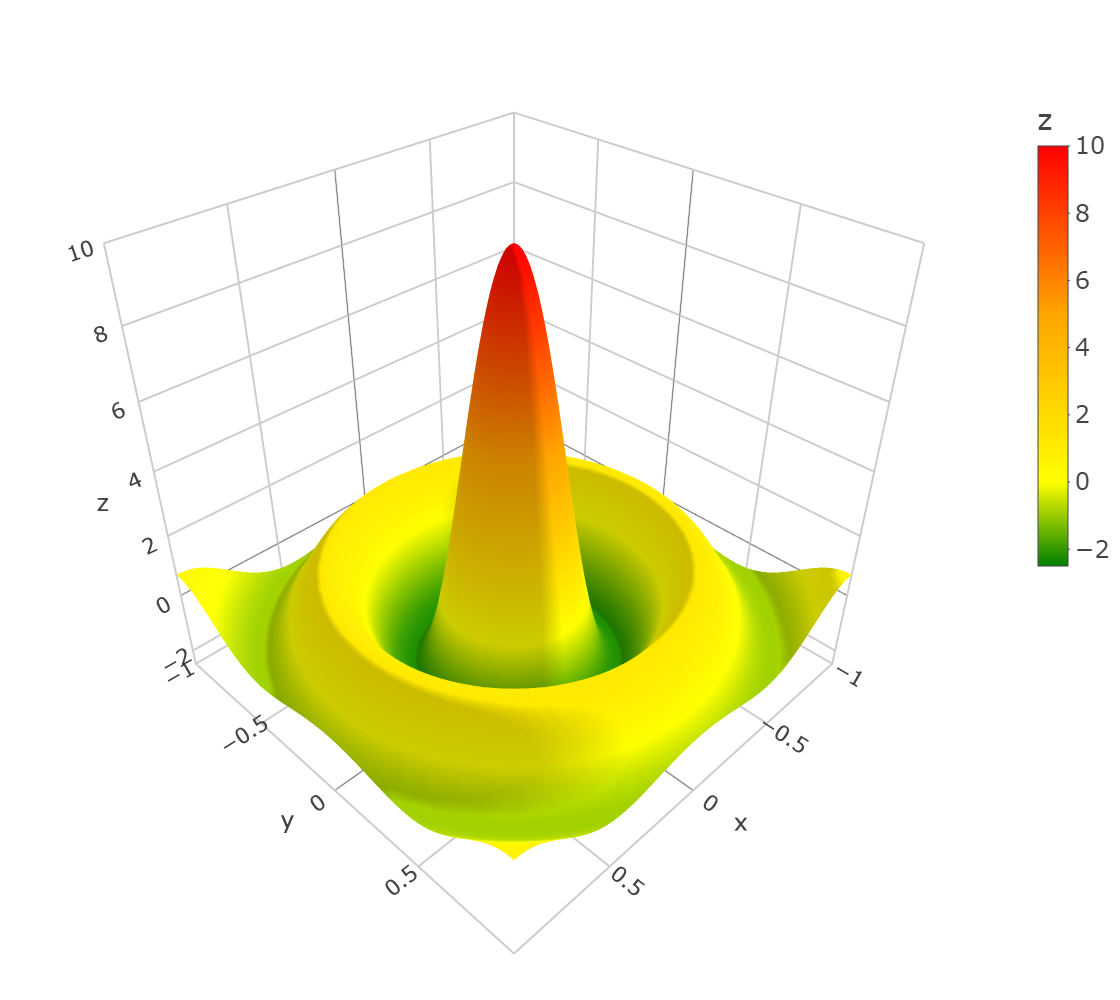}
    }
    \subfigure[A ULCV estimate $z=f^*(x,y)$]{
        \includegraphics[width=0.45\textwidth]{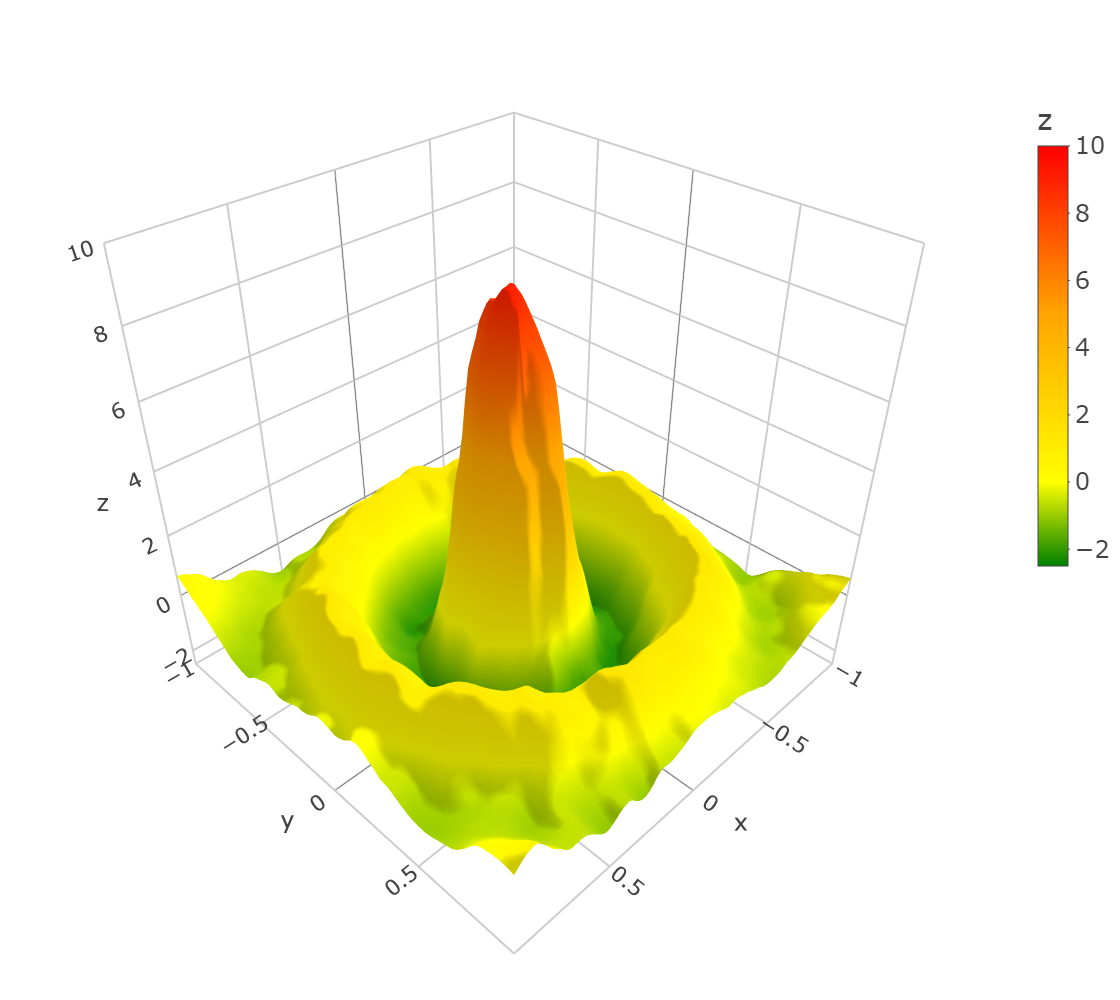}
    }
    \vspace{-12pt}
    \caption{The estimated function (left) and a result of ULCV estimator (right) in Example~2}\label{ex2fun}
\end{figure}

The i.i.d. design points were generated
with independent polar coordinates $(\rho,\varphi)$, where $\rho$ was drawn with the density proportional
to $r^2(2-r)^{1/10}$, $0\le r \le 2$, and $\varphi$ was uniformly distributed on $[0,2\pi)$. The distribution of the design points was restricted on
$\Theta = [-1,1]\times [-1,1]$, i.e., the design points that did not fall into $\Theta$ were excluded, keeping the total number of collected points equal to 5000, as in the other simulation examples.
One draw of the design points is depicted in Fig.~\ref{ex2mse}. The estimated function $f(x,y)$ and a computed ULCV estimate are depicted in Fig.~\ref{ex2fun}.

%
%
%

The results are presented in Fig.~\ref{ex2mse}. The ULCV estimator was the best among the three considered ones both for MSE and MaxE accuracy measures. In particular, the ULCV estimator was better than the NW one: MSE 0.2803 (0.2718,  0.2898) vs.
0.2870  (0.2774,  0.2974), $p<0.0001$; MaxE 2.505  (2.072, 3.140) vs.
2.695 (2.303, 3.361), $p<0.0001$.
In this example, the ULCM estimator had lower MaxE and higher MSE than the NW estimator did.

\bigskip

\newpage
\noindent\textbf{5.3. Example 3}

\begin{figure}[!h]
    \centering
    \subfigure[Design points] {
        \includegraphics[width=0.35\textwidth]{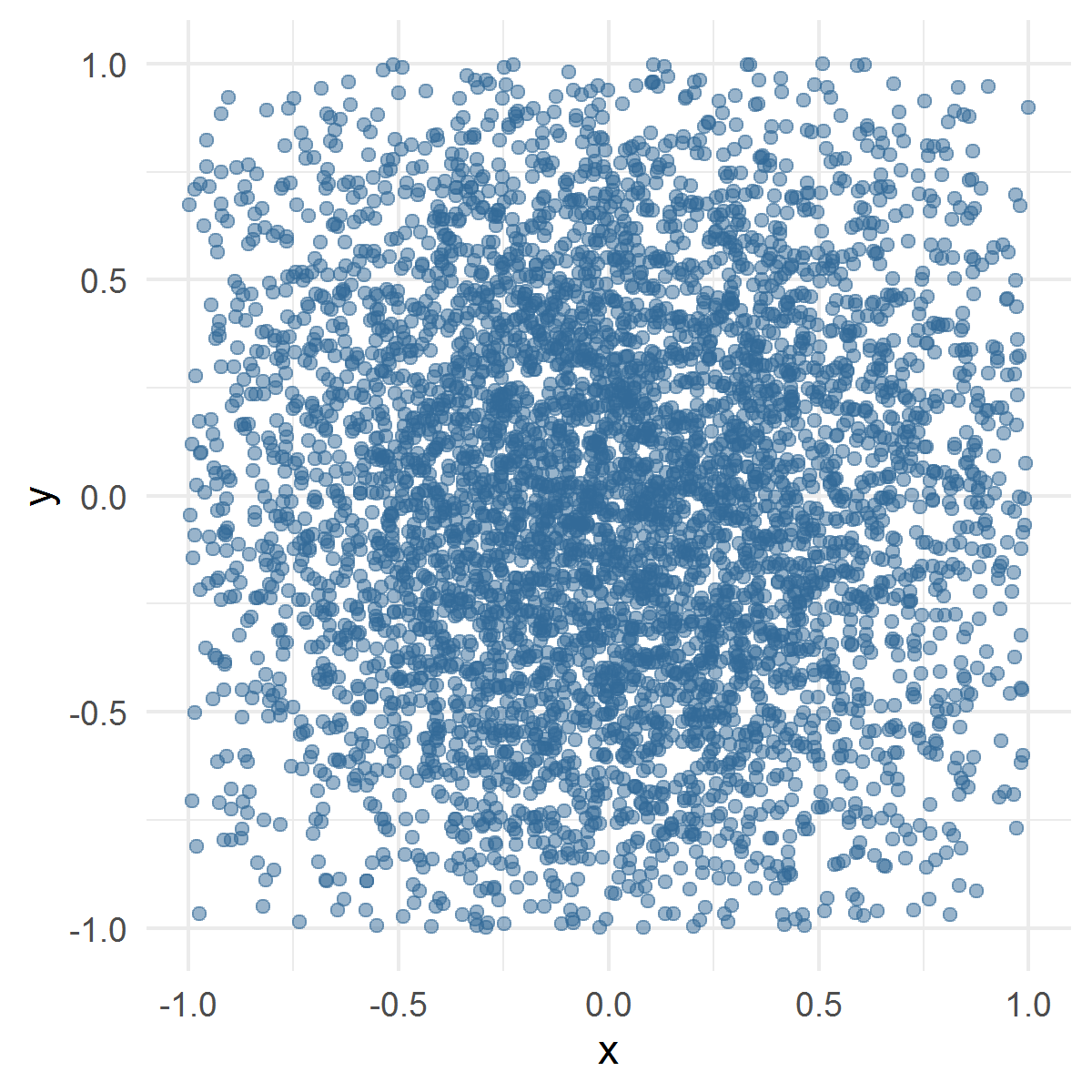}
    }
    \subfigure[MSE]{
        \includegraphics[width=0.25\textwidth]{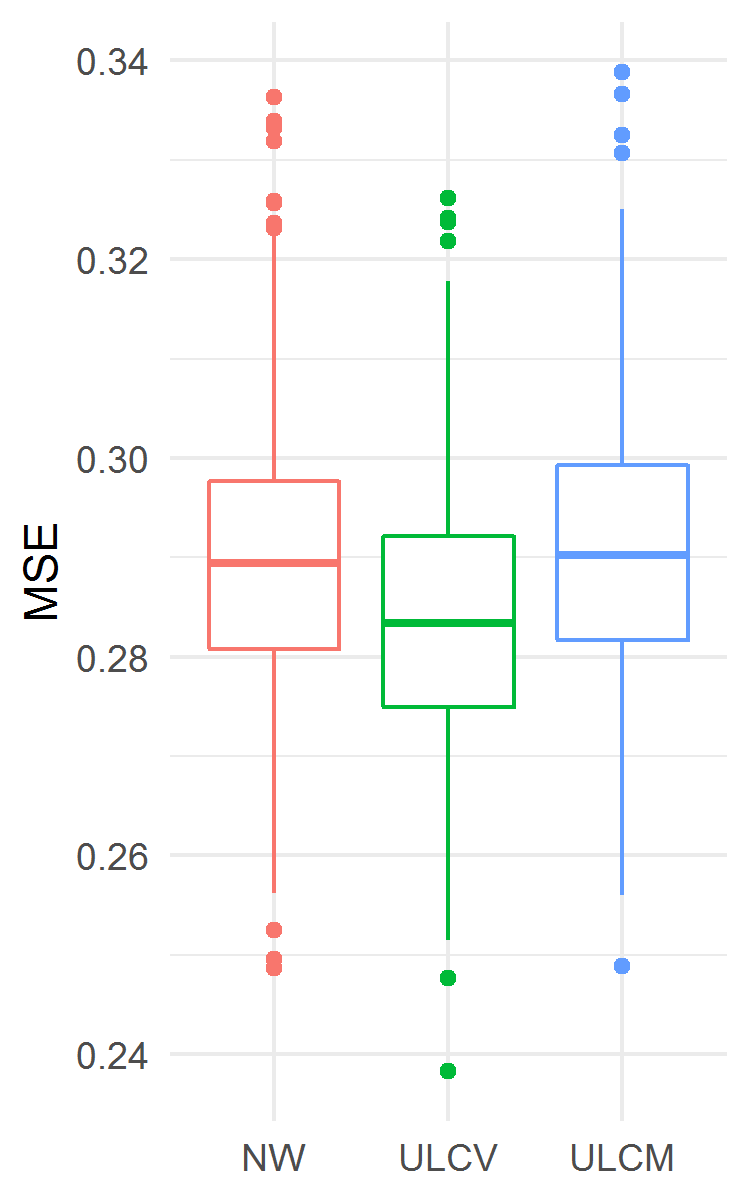}
    }
    \subfigure[MaxE]{
        \includegraphics[width=0.25\textwidth]{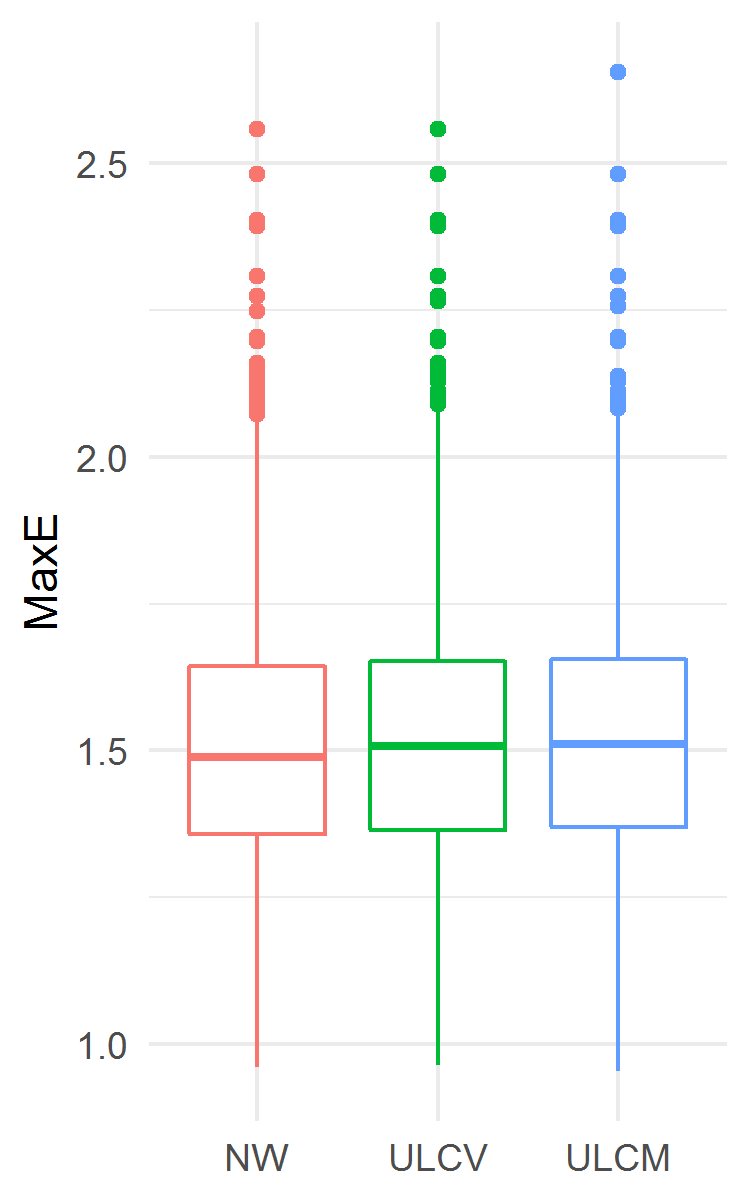}
    }
    \vspace{-12pt}
    \caption{Design points in an Example~3 experiment (left), the mean-square errors (middle), and maximal absolute errors (right) in Example~3}\label{ex3mse}
\end{figure}

In this example, we approximate the same nonrandom regression function
(\ref{funex2}) as in Example~2. The only difference of this example from Example~2 is that here
the coordinates of the design points were generated as independent normal random variables with mean 0 and standard deviation 1/2.
As above, the distribution of the design points was restricted on
$\Theta = [-1,1]\times [-1,1]$.
One draw of the design points is depicted in Fig.~\ref{ex3mse}.

%
%
%


The results are presented in Fig.~\ref{ex3mse}.
The ULCV estimator was the best one in terms of MSE, in particular, it was better than NW: MSE
0.2834  (0.2750, 0.2922) vs.
0.2895  (0.2808,  0.2977), $p<0.0001$.
But ULCV was worse than NW in terms of MaxE:
1.507  (1.364,   1.653) vs.
1.488 (1.357,   1.643), $p<0.0001$.
In this example, the ULCM estimator was the worst one both for MSE and MaxE.
However, from a practical point of view, the three estimators demonstrated similar accuracy in terms of MaxE.

\bigskip
\begin{center}
{\bf 6. Real data application}
\end{center}

\begin{figure}[!h]
    \centering
    \subfigure[Earthquakes] {
        \includegraphics[width=0.35\textwidth]{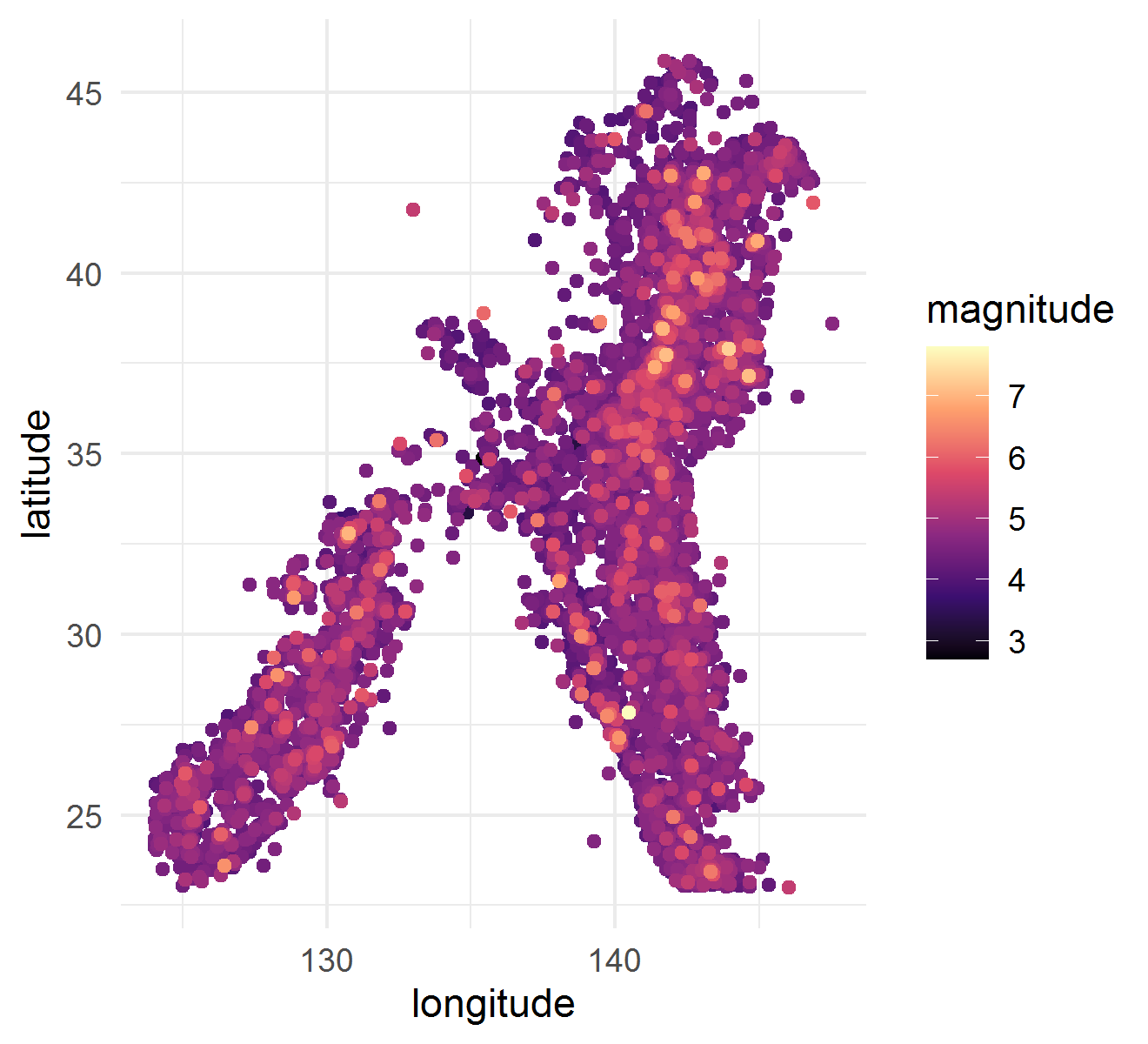}
    }
    \subfigure[MSE]{
        \includegraphics[width=0.25\textwidth]{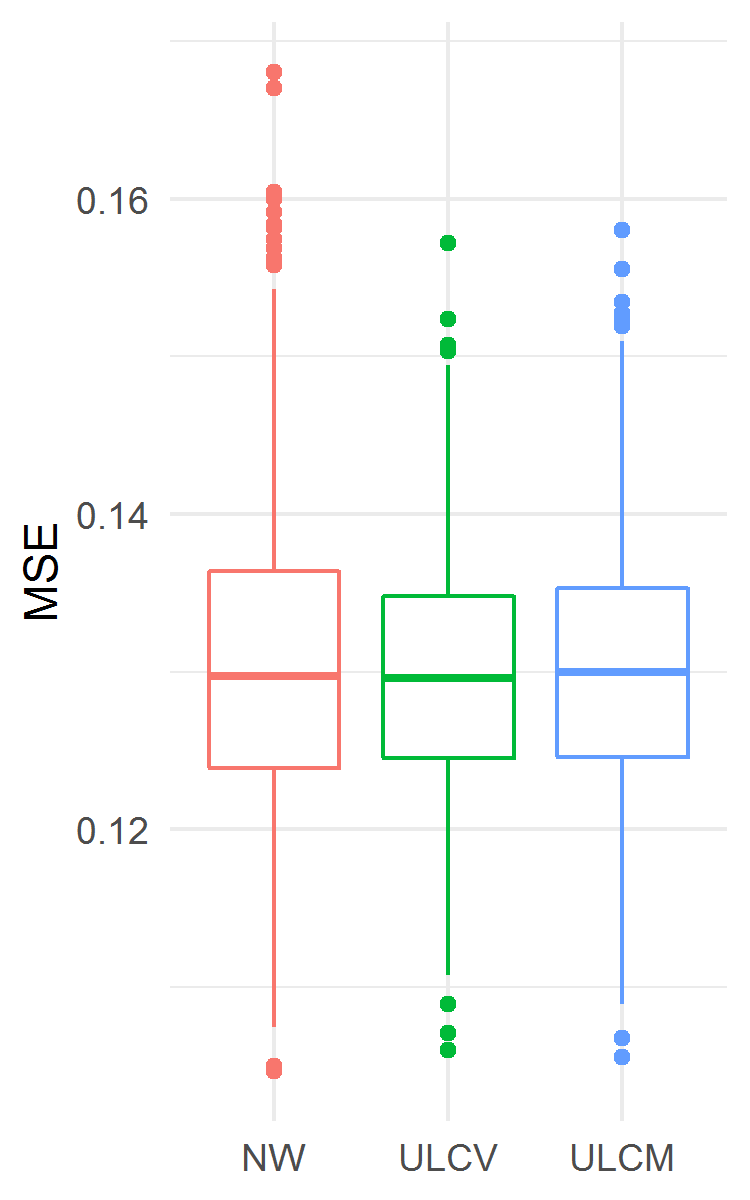}
    }
    \subfigure[MaxE]{
        \includegraphics[width=0.25\textwidth]{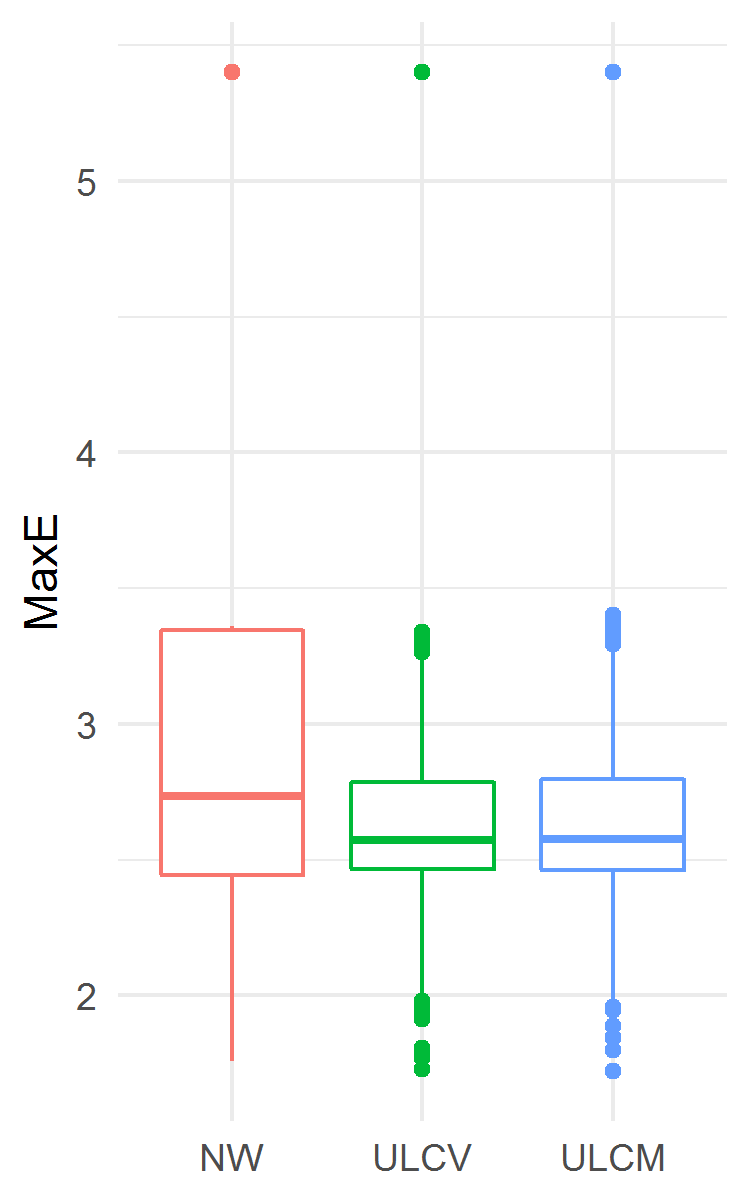}
    }
    \vspace{-12pt}
    \caption{Observed earthquakes events (left), the mean-square errors (middle), and maximal absolute errors (right)}\label{eqmse}
\end{figure}

\begin{figure}[!h]
    \centering
    \subfigure[The NW estimate $z=\hat f(x,y)$]{
        \includegraphics[width=0.45\textwidth]{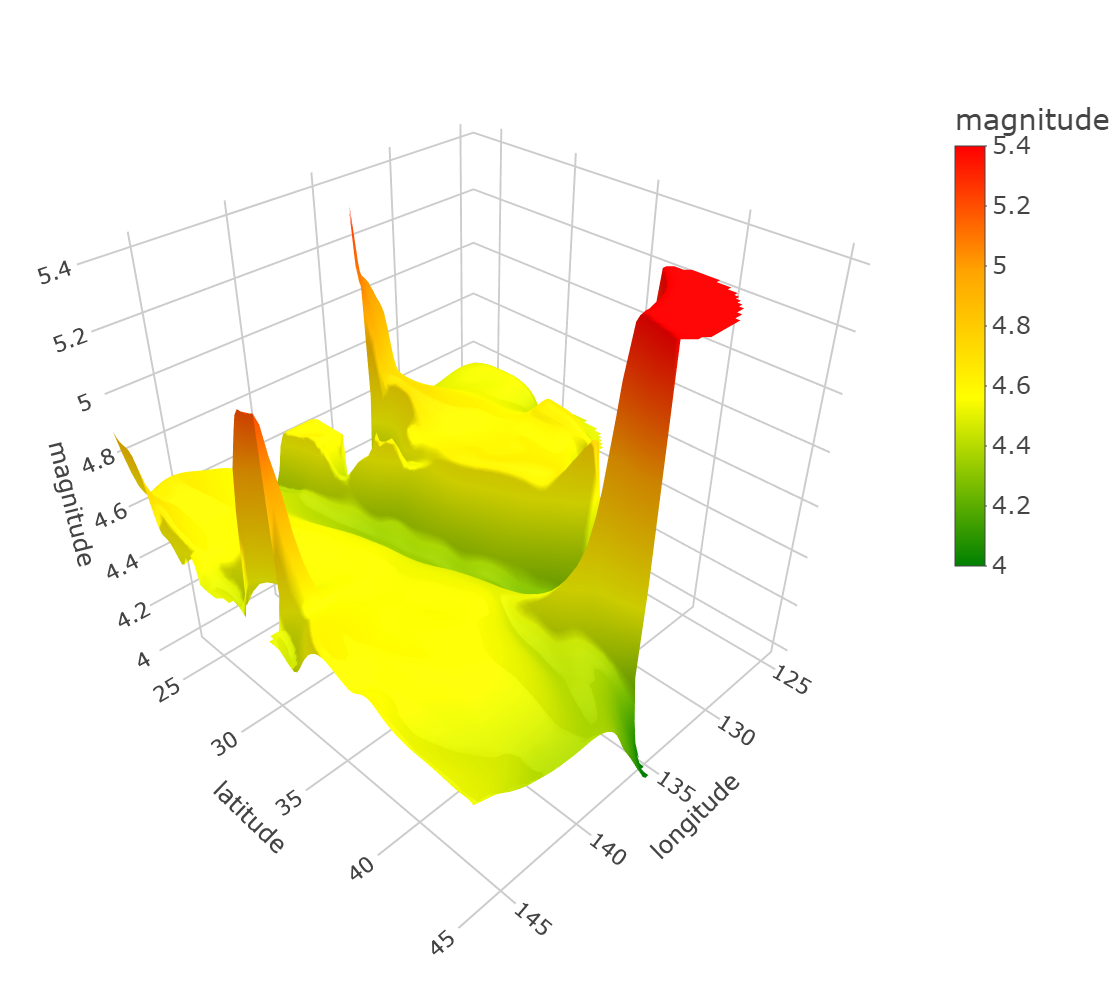}
    }
    \subfigure[The ULCV estimate $z=f^*(x,y)$]{
        \includegraphics[width=0.45\textwidth]{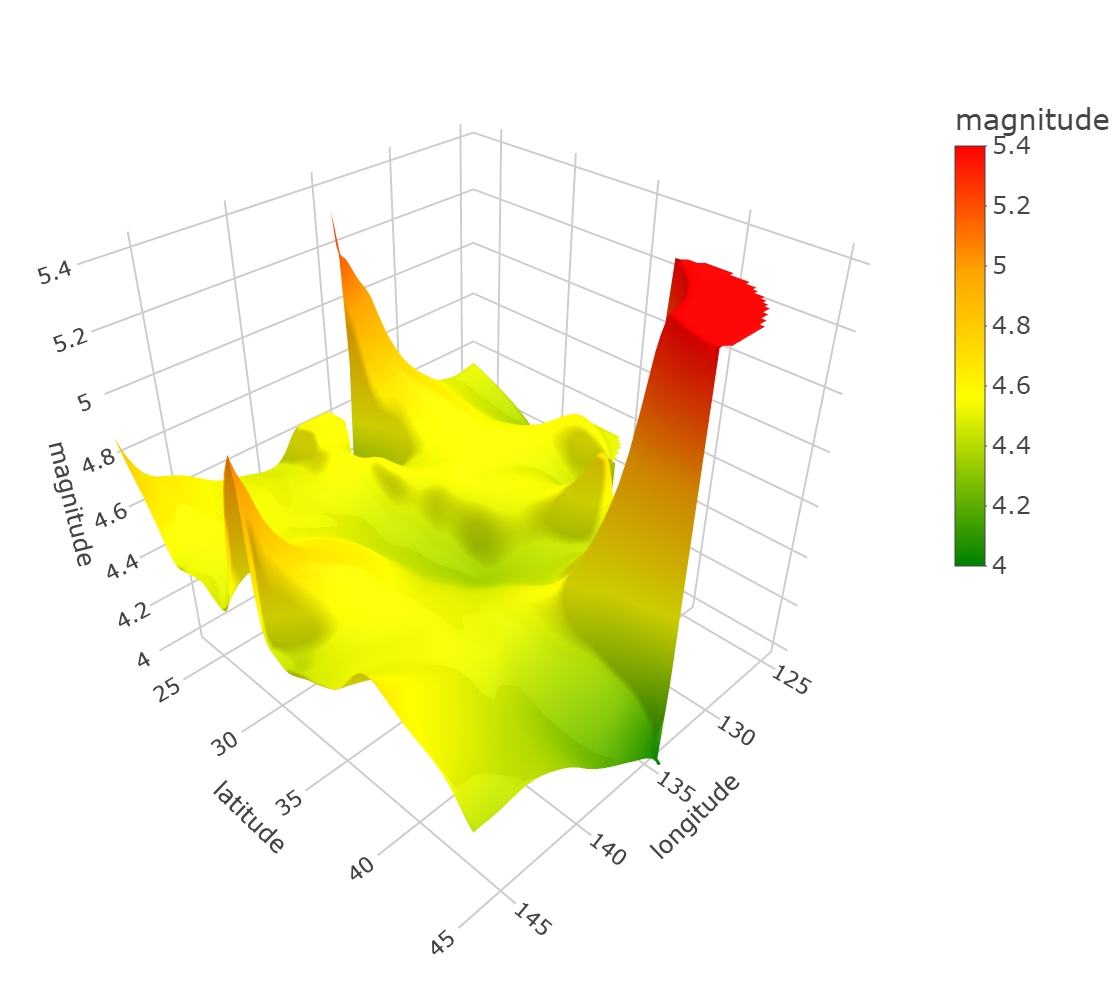}
    }
    \vspace{-12pt}
    \caption{The estimated mean magnitude for the NW (left) and ULCV (right) estimators}  \label{eq7est}
\end{figure}

In this section, we compared the new ULCV and ULCM estimators with the NW one in the application to the data on earthquakes in Japan that happened in 2012--2021 (data retrieved from ANSS Comprehensive Earthquake Catalog, 2022). Each of the 10184 collected earthquake events was described by its coordinates (longitude and latitude) and its magnitude (ranging from 2.7 to 7.8).
The collected events are presented in Fig.~\ref{eqmse}.
The goal of the application of the estimators was to accurately estimate the mean magnitude depending on the coordinates.
As in the simulation examples above, we did 1000 runs, in each of which the data were randomly divided into the training (80\%) and validation (20\%) sets. For each of the tested algorithms, on the training set, the optimal $\varepsilon$ was calculated by 10-fold cross-validation minimizing the average of mean-square errors.
The $\varepsilon$ was selected from 20 values located on the logarithmic grid from 1 to 10.
The random partitioning for the cross-validation was the same for all the tested algorithms.
The difference of the computations of this section with those of Sec.~5 was that we did not know the true value of the estimated function, therefore, we had to estimate the maximal error (MaxE) on the validation set in each run, not on true values of the estimated function.
Besides, since the domain of the coordinates of the events is nonrectangular while the epmloyed domain partitioning algorithms (Voronoi cells algorithm and coordinate-wise medians algorithm) calculated the squares of the cells for a rectangular domain, we bounded the squares of the cells from above by 1 in order to avoid overweighting of the corresponding observations. The resulting estimates of the NW and ULCV estimators are depicted in Fig.~\ref{eq7est}, where, for each estimator, the value of $\varepsilon$ was chosen as the median of those chosen in the 1000 runs.

The results are presented in Fig.~\ref{eqmse}. The ULCV estimator was the best among the three considered ones both for MSE and MaxE accuracy measures. In particular, the ULCV estimator was better than the NW one: MSE
 0.1296  (0.1245,    0.1348) vs.
0.1297   (0.1239,    0.1364),  $p<0.0001$;
MaxE
2.573   (2.464,   2.785) vs.
2.736   (2.442,    3.346), $p=0.0005$.
In this example, the ULCM estimator yielded lower MaxE and higher MSE than the NW estimator did.
However, from a practical point of view, the three estimators displayed similar median MSE.


\begin{center}
{\bf 7. Conclusion}
\end{center}

In this paper, for a wide class of nonparametric regression models with a multivariate random design,
 universal uniformly consistent kernel estimators are  proposed for the unknown random regression functions (random fields) of the corresponding multivariate argument. These estimators belong
to the class of local constant kernel estimators.
But in contrast to the vast majority of previously known results, traditional correlation conditions of design elements
are not needed for the consistency of the new estimators.
The design can be either fixed and not necessarily regular, or random and not necessarily
consisting of independent or weakly dependent random variables. With regard to design elements, the only condition  that is required is the dense filling of the regression function domain with the design points.

Explicit upper bounds are found for the rate of uniform convergence in probability of the new estimators to an unknown random regression function.
The only characteristic explicitly included in these estimators is the maximum diameter of the cells of
partition generated by the design elements,
and only convergence to zero in probability is required for the characteristic.
The advantage of this condition over the classical ones
is that it is insensitive to the forms of dependence of the design observations.
Note that this condition is, in fact, necessary, since only when the design densely fills the regression function domain, it is possible to reconstruct the regression function with a certain accuracy. As a corollary of the main result, we obtain a consistent estimator for the mean function of a continuous random process.

In the simulation examples of Section~5, the new estimators were compared
with Nadaraya--Watson estimators. In some of the examples,
the new estimators proved to be most accurate.
In Section~6, as an application of the new estimators, we studied the real data on the magnitudes of earthquakes in Japan, and
the accuracy of the new estimators was comparable to that of the Nadaraya-Watson ones.

\begin{center}
{\bf 8. Proofs}
\end{center}

{\it Proof of Theorem $1$.} Taking the relation (\ref{f2}) into account, one can obtain the following identity:
$$f^*_{n,\varepsilon}(t)=R_{n,\varepsilon}(f,t)+\nu_{n,\varepsilon}(t),
$$
where
$$R_{n,\varepsilon}(f,t):=J^{-1}_{n,\varepsilon}(t)\sum\limits_{i=1}^nf(X_{i})K_{\varepsilon}\left(t-X_{i}\right)
\Lambda(\Delta_i),$$
$$\nu_{n,\varepsilon}(t):=J^{-1}_{n,\varepsilon}(t)\sum\limits_{i=1}^nK_{\varepsilon}\left(t-X_{i}\right)\xi_{i}\Lambda(\Delta_i),$$
$$J_{n,\varepsilon}(t):=\sum\limits_{i=1}^nK_{\varepsilon}(t-X_{i})\Lambda(\Delta_i).
$$

Notice that, by virtue of the properties of $K$,
 the range of summation in the three sums above is equal to $\{i:\,\|t-X_{i}\|\le\varepsilon\}$.
This is the principal argument in the calculations below.

{\bf Lemma 1}. {\it The following estimate is valid:
\begin{eqnarray}\label{approx}
\sup_{t\in \Theta}|R_{n,\varepsilon}(f,t)-f(t)|\le\omega_f(\varepsilon).
\end{eqnarray}

The proof} is immediate from the identity
$$R_{n,\varepsilon}(f,t)=f(t)+J^{-1}_{n,\varepsilon}(t)\sum\limits_{i:\,\|t-X_{i}\|\le\varepsilon}(f(X_{i})-f(t))
K_{\varepsilon}\left(t-X_{i}\right)\Lambda(\Delta_i).
$$
$\hfill\Box$

{\bf Lemma 2}. {\it If $\delta_n\le\varepsilon\le\varepsilon_0$ then, for any $t\in \Theta$, the following relation is valid$:$
\begin{eqnarray*}\label{riemann}
J_{n,\varepsilon}(t)\ge \rho - k 2^k L  \delta_n \varepsilon^{-1}.
\end{eqnarray*}

Proof.}
We have
$$J_{n,\varepsilon}(t)=\int  g(y) dy,$$
where
$$g(y)=\sum\limits_{i=1}^nK_{\varepsilon}(t-X_{i})
\,I(y\in \Delta_i);$$
here $I(\cdot)$ is the indicator of an event.
Since
\begin{equation}\label{jne}
|K_{\varepsilon}\left(t-x\right)-K_{\varepsilon}\left(t-y\right)|\le kL\varepsilon^{-k-1}\delta_n,
 \  \mbox{ whenever }\  \|x-y\|\le \delta_n,
\end{equation}
we have
$g(y)\ge
K_{\varepsilon}(t-y)-kL\varepsilon^{-k-1}\delta_n
$ for all $y\in\Theta$.
Hence, by (\ref{rho}),
\begin{eqnarray*}
J_{n,\varepsilon}(t)&\ge&
\int_{y\in\Theta:\ \|t-y\|\le \varepsilon}
\big(K_{\varepsilon}(t-y)-kL\varepsilon^{-k-1}\delta_n\big) dy
\\
&\ge&
\rho - k 2^kL  \delta_n \varepsilon^{-1}.
\end{eqnarray*}

The lemma is proved.$\hfill\Box$

{\bf Lemma 3}. {\it For every $y>0$ and $\varepsilon\in(0, \varepsilon_0]$, on the subset of elementary events defined by the relation  $\delta_n/\varepsilon\le
\min\{1,\, \rho(k2^{k+1}L)^{-1}\}$, the following upper bound is valid$:$
\begin{eqnarray}\label{noise}
{\bf P}_{\cal F}\left (\Big(\delta_n^{k/2}\varepsilon^{-k((1/2)+(1/p))}\Big)^{-1} \sup\limits_{t\in \Theta}|\nu_{n,\varepsilon}(t)|>y\right )\le
G(k,p)\, \rho^{-p}\, M_p\, L^{p/2} y^{-p},
\end{eqnarray}
where the symbol ${\bf P}_{\cal F}$ denotes the conditional probability given the $\sigma$-field ${\cal F}$ generated by the design $\{X_i\}$ and the paths of the random field $f(\cdot)$$;$ here
$$G(k,p)< (p-1)^{p/2} 2^{p(k+(3/2))}
\,
\left( 1
+\frac{k}{2^{(p-k)/(p+1)} - 1}
\right)^{p+1}.
$$
}

{\it Proof.} Under the condition $\delta_n/\varepsilon\le \rho(k2^{k+1}L)^{-1}$,  by virtue of Lemma 2 the simple inequality
$|\nu_{n,\varepsilon}(t)|\le 2\rho^{-1}|\mu_{n,\varepsilon}(t)| $  is valid, where
$$
\mu_{n,\varepsilon}(t) := \sum\limits_{i=1}^nK_{\varepsilon}(t-X_{i})\Lambda(\Delta_i)\xi_{i}.
$$
The distribution tail of
$\sup_{t\in \Theta}|\mu_{n,\varepsilon}(t)|$
will be estimated by   Kolmogorov's  {\it dyadic chaining} 
 which has been used  to estimate the tail probability of the sup-norm of a stochastic processes
having continuous paths with probability 1.

Without loss of generality we will assume that $\Theta\subset [0,1]^k$.
We first note that
the set $\Theta$ under the supremum sign above can be replaced with the subset of dyadic rational points
${\cal R}=\cup_{l\ge 1} {\cal R}_l$,
where
$${\cal R}_l=\{(j_1/2^l,\dots,j_k/2^l):\  j_1=1,\dots, 2^l-1;\dots;
j_k=1,\ldots, 2^l-1 \}.$$
Thus,
$$\sup_{t\in \Theta}|\mu_{n,\varepsilon}(t)|
\le \sup_{t\in {\cal R}}|\mu_{n,\varepsilon}(t)|$$
$$\le \max_{t\in {\cal R}_m}|\mu_{n,\varepsilon}(t)| +
\sum_{l=m+1}^\infty
\sum_{r=1}^k
\max_{t\in {\cal R}_l}
\big|\mu_{n,\varepsilon}(t+2^{-l}e_r)-\mu_{n,\varepsilon}(t)\big|,
$$
where $m$ is some natural number that will be chosen later, and $e_r$ is the
$k$-dimensional vector with the $r$-th component 1 and other components 0.

Hence,
\begin{eqnarray}
{\bf P}_{\cal F}(\sup_{t\in \Theta}|\mu_{n,\varepsilon}(t)|&>&y)
\notag \\
&\le&
{\bf P}_{\cal F}(\max_{t\in {\cal R}_m}|\mu_{n,\varepsilon}(t)|
>a_m y)
\notag \\
 &&+
\sum_{l=m+1}^\infty
\sum_{r=1}^k
{\bf P}_{\cal F}\Big(
\max_{t\in {\cal R}_l}
\big|\mu_{n,\varepsilon}(t+2^{-l}e_r)-\mu_{n,\varepsilon}(t)\big|
>a_l y/k\Big)
\notag \\
\label{req1}
&\le&
\sum_{t\in {\cal R}_m}
{\bf P}_{\cal F}(|\mu_{n,\varepsilon}(t)|
>a_m y)
\\
&&+
\sum_{l=m+1}^\infty
\sum_{r=1}^k
\sum_{t\in {\cal R}_l}
{\bf P}_{\cal F}\Big(
\big|\mu_{n,\varepsilon}(t+2^{-l}e_r)-\mu_{n,\varepsilon}(t)\big|
>a_l y/k\Big)
,
\notag
\end{eqnarray}
where $a_m,a_{m+1},\dots$ is a sequence of positive numbers such that
$a_m+a_{m+1}+\cdots=1$.

In order to estimate the probability
${\bf P}_{\cal F}(|\mu_{n,\varepsilon}(t)| >a_m y)$,
we use Rio's  martingale inequality
(Rio 2009, Theorem 2.1)
\begin{eqnarray}\label{rosen}
{\bf E}\Big| \sum^n_{i=1}\eta_i \Big|^p\le \left(
(p-1)
\sum^n_{i=1}\big({\bf E}|\eta_i|^p\big)^{2/p}\right)^{p/2},
\end{eqnarray}
where $\{\eta_i\}$ is a martingale-difference sequence with
finite moments of order $p\ge 2$.

Now, put
$$\eta_i:=K_{\varepsilon}(t-X_{i})\Lambda(\Delta_i)\xi_{i}.$$
From
(\ref{rosen}) and the simple
upper bounds
$$K_{\varepsilon}(t-X_{i})\Lambda(\Delta_i)\le
\frac{L}{\varepsilon^k}\delta^k_n,$$
$$\sum_i K_{\varepsilon}(t-X_{i})\Lambda(\Delta_i)\le
\frac{L}{\varepsilon^k}(2\varepsilon+2\delta_n)^k
$$
we then obtain that, with probability $1$,
$$\sum^n_{i=1}\big({\bf E}_{\cal F}|\eta_i|^p\big )^{2/p}\le
M_p^{2/p}\, 2^k L\,  (1+\delta_n/\varepsilon)^{k}
\,(\delta_n/\varepsilon)^{k}.
$$

Under the restriction $\delta_n\le\varepsilon\le 1$, the last inequality and (\ref{rosen}) imply that
\begin{eqnarray}\label{final1}
{\bf P}_{\cal F}(|\mu_{n,\varepsilon}(t)| >a_m y)
\le \frac{{\bf E}_{\cal F} (|\mu_{n,\varepsilon}(t)|^p)}{(a_m y)^p}
\le G_1
\frac{(\delta_n/\varepsilon)^{kp/2}}{(a_m y)^p} \quad\mbox{a.s.},
\end{eqnarray}
where
$$G_1 = (p-1)^{p/2}\, 2^{kp} L^{p/2}\,  M_p  .$$

Now let us estimate
${\bf P}_{\cal F}\Big(
\big|\mu_{n,\varepsilon}(t+2^{-l}e_r)-\mu_{n,\varepsilon}(t)\big|
>a_l y/k\Big)$. In order to do it we use (\ref{rosen}) with
$$\eta_i:= \Big(K_{\varepsilon}(t-X_{i}+2^{-l}e_r)-
K_{\varepsilon}(t-X_{i}) \Big)
\Lambda(\Delta_i)\xi_{i}.$$
We have
$$ \Big|K_{\varepsilon}(t-X_{i}+2^{-l}e_r)-
K_{\varepsilon}(t-X_{i}) \Big| \Lambda(\Delta_i) \le
\frac{L}{\varepsilon^{k+1}}\ 2^{-l}\ \delta_n^k,
$$
$$
\sum_i \Big|K_{\varepsilon}(t-X_{i}+2^{-l}e_r)-
K_{\varepsilon}(t-X_{i}) \Big|\Lambda(\Delta_i)\le
\frac{L}{\varepsilon^{k+1}}\ 2^{-l+1}\ (2\varepsilon+2\delta_n)^k.
$$
Thus
$$\sum^n_{i=1}\big({\bf E}_{\cal F}|\eta_i|^p\big )^{2/p}\le
M_p^{2/p}\, 2^k L\, 2^{-2l+1}\, (1+\delta_n/\varepsilon)^{k}
\, (\delta_n/\varepsilon)^{k}.
$$
Again, under the restriction $\delta_n\le\varepsilon\le 1$, the last inequality and (\ref{rosen}) imply
\begin{eqnarray}\label{final2}
 {\bf P}_{\cal F}\Big(
\big|\mu_{n,\varepsilon}(t+2^{-l}e_r)-\mu_{n,\varepsilon}(t)\big|
>a_l y/k\Big)
\le \frac{G_2}{k}
\ \frac{(\delta_n/\varepsilon)^{kp/2}\ \varepsilon^{-p} \ 2^{-lp}}{(a_l y)^p},
\end{eqnarray}
where
$$G_2= 2^{p/2}\,k^{p+1}\, (p-1)^{p/2}\, 2^{kp} L^{p/2}\,  M_p =2^{p/2}\,k^{p+1}\, G_1.$$

Combining (\ref{req1}), (\ref{final1}), and (\ref{final2}), we obtain
$${\bf P}_{\cal F}(\sup_{t\in \Theta}|\mu_{n,\varepsilon}(t)|>y)$$
$$< y^{-p} (\delta_n/\varepsilon)^{kp/2}
\left( G_1 2^{km} a_m^{-p}
+G_2\, \varepsilon^{-p} \sum_{l=m+1}^\infty 2^{-(p-k)l}\, a_k^{-p}
\right).
$$
The optimal sequence $a_l$ minimizing the right-hand side of this inequality
is as follows: $a_m=c (G_1 2^{km})^{1/(p+1)} $ and
$a_l=c \big(G_2\, \varepsilon^{-p} \, 2^{-(p-k)l}\big)^{1/(p+1)}$  for $l=m+1,m+2,\dots$, where the coefficient $c$ is defined by the relation
$a_m+a_{m+1}+\cdots=1$.
For this sequence, we get
$${\bf P}_{\cal F}(\sup_{t\in \Theta}|\mu_{n,\varepsilon}(t)|>y)$$
$$< y^{-p} (\delta_n/\varepsilon)^{kp/2}
\left( (G_1 2^{km})^{1/(p+1)}
+\sum_{l=m+1}^\infty \Big(G_2\, \varepsilon^{-p}  2^{-(p-k)l}
\Big)^{1/(p+1)}
\right)^{p+1}.
$$
Now, put
$m = \lceil- \log_2 \varepsilon\rceil$, where $\lceil a \rceil$ is the minimal
integer greater than or equal to $a$. Then
$${\bf P}_{\cal F}(\sup_{t\in \Theta}|\mu_{n,\varepsilon}(t)|>y)$$
$$< y^{-p}\, \delta_n^{kp/2}\, \varepsilon^{-k(1+ (p/2))}
\left( (2 G_1 )^{1/(p+1)}
+\big(G_2  2^{-(p-k)} \big)^{1/(p+1)} \sum_{l=0}^\infty 2^{-(p-k)l/(p+1)}
\right)^{p+1}
$$
$$< y^{-p}\, \delta_n^{kp/2}\,\varepsilon^{-k(1+ (p/2))}
2^{p/2} \,G_1
\left( 1
+\frac{k}{2^{(p-k)/(p+1)} - 1}
\right)^{p+1}.
$$
This yields the statement of the lemma with
$$G(k,p)=2^p\,2^{p/2} \,\frac{G_1}{L^{p/2}\,  M_p}
\left( 1
+\frac{k}{2^{(p-k)/(p+1)} - 1}
\right)^{p+1}.
$$
Lemma 3 is proved.$\hfill\Box$

The statement of Theorem 1 follows from Lemmas 1--3.

{\it Proof of Theorem} 2. First of all, notice that condition (\ref{supmoment}) and  Lebesgue's dominated convergence theorem imply
the relation
\begin{equation}\label{LLN1}
\lim_{\nu\to 0}{\bf E}\omega_f(\nu)=0.
\end{equation}
It is clear that the relation (\ref{LLN1}) implies  the uniform law of large numbers  for
independent copies of the a.s. continuous random process $f(t)$, i.e.,
$$\sup\limits_{t\in [0,1]^k}|\overline{f}_{N}(t)-{\bf E}f(t)|\stackrel{p}\to 0
$$
as $N\to \infty$, where $\overline f_{N}(t):=\frac{1}{N}\sum\limits_{j=1}^Nf_j(t)$.
Put
\begin{equation}\label{resudial}
\Delta_{n,\varepsilon,j}:=\sup_{t\in [0,1]^k}| f^*_{n,\varepsilon,j}(t)-f_j(t)|.
\end{equation}
So, to prove (\ref{LLN}) we need only  to verify the following version of the law of large numbers
for independent copies of the residuals defined in (\ref{resudial}):
$$
\frac{1}{N}\sum_{j=1}^N\Delta_{n,\varepsilon,j}\stackrel{p}\to 0,
$$
but only for the sequences $\varepsilon$ and $N$ chosen in (\ref{Neps}).

Introduce the following events:
$$A_{n,\varepsilon,j}:=\{\delta_{n,j}\le\varepsilon\min\{1,\, \rho(k2^{k+1} L)^{-1}\}\},\,\,\,j=1,\ldots,N,
$$
where the sequence $\varepsilon\equiv \varepsilon_n\to 0$ meets (\ref{Neps}). (It is evident that such a sequence exists.)
For any positive $\nu$ we have
\begin{equation}\label{ineq}
{\bf P}\left \{\frac{1}{N}\sum_{j=1}^N\Delta_{n,\varepsilon,j}>\nu\right\}\le {\bf P}\left \{\frac{1}{N}
\sum_{j=1}^N\Delta_{n,\varepsilon,j}I(A_{n,\varepsilon,j})>\nu\right\}+N{\bf P}(\overline{A_{n,\varepsilon,1}}).
\end{equation}
Next, from Theorem 1 we obtain
$${\bf E}\Delta_{n,\varepsilon,j}I(A_{n,\varepsilon,j})\le {\bf E}\omega_f(\varepsilon)+
\int\limits_0^{\infty}{\bf P}\left( \zeta_n(\varepsilon)>y,\, \delta_n\le\varepsilon\min\{1,\, \rho(k2^{k+1} L)^{-1}\}\right)dy
$$
$$\le {\bf E}\omega_f(\varepsilon)+\gamma_n +
\int\limits_{\gamma_n}^{\infty}{\bf P}\left( \zeta_n(\varepsilon)>y,\, \delta_n\le\varepsilon\min\{1,\, \rho(k2^{k+1} L)^{-1}\}\right)dy
$$
$$\le {\bf E}\omega_f(\varepsilon)+\tilde C \gamma_n,
$$
where $\tilde C:=1+G(k,p)\, \rho^{-p}\, M_p\, L^{p/2}$ and $\gamma_n:=\left(\varepsilon^{-k(p/2+1)}\,{\bf E}\delta_n^{kp/2}\right)^{1/p}$. It remains to apply Markov's inequality for the first probability on the right-hand side of (\ref{ineq})
and use the limit relations (\ref{Neps}) and the last estimate. Theorem 2 is proved.$\hfill\Box$


\bigskip

\begin{center}
{\bf Acknowledgments}
\end{center}
The authors are deeply grateful to Professor I.A. Ibragimov for his useful remarks.
In addition, the authors thank the anonymous referee whose comments contributed to a better presentation of this study.

{\bf Data availability statement.}
The data required to reproduce the above findings are available to download from
https://earthquake.usgs.gov/data/comcat/
(ANSS Comprehensive Earthquake Catalog, 2022).

\begin{center}
{\bf References}
\end{center}

\textcolor{blue}{
Ahmad, I. A. and Lin, P.-E. (1984), `Fitting a multiple regression function', \emph{J. Statist. Plann. Infer.} 9, 163--176.
}

\textcolor{blue}{
ANSS Comprehensive Earthquake Catalog, 2022. In:
U.S. Geological Survey, Earthquake Hazards Program, 2017, Advanced National Seismic System (ANSS) Comprehensive Catalog of Earthquake Events and Products: Various, https://doi.org/10.5066/F7MS3QZH.
Data retrieved September 4, 2022 from https://earthquake.usgs.gov/data/comcat/
}

Benhenni, K., Hedli-Griche, S.,  and Rachdi, M. (2010),
`Estimation of the regression operator from functional fixed-design with correlated errors',
\emph{J.  Multivar. Anal.} 101, 476--490.

\textcolor{blue}{
Benelmadani, D., Benhenni, K.,  and  Louhichi, S. (2020),  `Trapezoidal rule and sampling designs for the nonparametric estimation of the regression function in models with correlated errors', \emph{Statistics} 54, 59--96.
}

Borisov, I.S., Linke, Yu.Yu., and Ruzankin P.S. (2021), `Universal weighted kernel-type estimators for some class of regression models',
{\it Metrika} 84, 141--166.

Brown, L.D. and Levine, M. (2007), `Variance estimation in nonparametric regression via the difference sequence method',
\emph{Ann. Statist.} 35, 2219--2232.

\textcolor{blue}{
Chan, N. and   Wang, Q. (2014), `Uniform convergence for Nadaraya-Watson estimators with nonstationary data', \emph{Econometric Theory} 30, 1110--1133.
}

\textcolor{blue}{
Chu, C. K. and  Deng, W.-S. (2003),   `An interpolation method for adapting to sparse design in multivariate nonparametric regression', \emph{J. Statist. Plann. Inference} 116,  91--111.
}

Einmahl, U. and Mason, D.M. (2005), `Uniform in bandwidth consistency of kernel-type
function estimators', \emph{Ann. Statist.}  33, 1380--1403.

Fan, J. and Gijbels, I.  (1996), \emph{Local Polynomial Modelling and its Applications},
London: Chapman and Hall.

\textcolor{blue}{
Fan, J. and  Yao, Q. (2003), \emph{Nonlinear time series nonparametric and parametric methods},  Springer.
}

\textcolor{blue}{
Gao, J.,  Kanaya, S.,  Li, D., and  Tjostheim, D. (2015), `Uniform consistency for nonparametric estimators in null recurrent time series', \emph{Econometric Theory} 31, 911--952.
}

Gasser, T. and Engel, J. (1990), `The choice of weghts in kernel regression estimation', \emph{Biometrica} 77, 277-381.

\textcolor{blue}{
Georgiev, A. A. (1988), `Consistent nonparametric multiple regression: The fixed design case', \emph{J. Multivariate Anal.} 25, 100--110.
}

\textcolor{blue}{
Georgiev, A. A. (1990), `Nonparametric multiple function fitting', \emph{Stat.  Probab. Lett.} 10, 203--211.
}

\textcolor{blue}{
Georgiev, A. A. (1989),  `Asymptotic properties of the multivariate Nadaraya-Watson regression function estimate: The fixed design case', \emph{Stat.  Probab. Lett.} 7, 35--40.
}

\textcolor{blue}{
Gu, W., Roussas, G. G., and  Tran, L. T. (2007),  `On the convergence rate of fixed design regression estimators for negatively associated random variables', \emph{Stat.  Probab. Lett.} 77, 1214--1224.
}

 Gy\"orfi, L.,  Kohler, M.,
 Krzyzak, A.,  and Walk, H. (2002),
\emph{A Distribution-Free
Theory of Nonparametric
Regression},  New York: Springer.

Hall, P. and Heyde, C. C., (1980),
{\it Martingale limit theory and its application.}
Academic Press.

Hall, P., M\"uller, H.-G.,  and Wang, J.-L. (2006), `Properties of principal component
methods for functional and longitudinal data analysis', \emph{Ann. Statist.} 34, 1493--1517.

Hansen, B.E. (2008),
`Uniform convergence rates for kernel estimation with dependent data',
\emph{Econometric Theory}  24,  726--748.

H\"ardle, W. (1990), \emph{Applied Nonparametric Regression}, New York: Cambridge
University Press.

\textcolor{blue}{
He, Q. (2019),  `Consistency of the Priestley--Chao estimator in nonparametric regression model with widely orthant  dependent errors', \emph{J. Inequal. Appl.} 64, 2--13.
}

\textcolor{blue}{
Hsing, T. and    Eubank, R.  (2015), \emph{Theoretical foundations of functional data
analysis, with an introduction to linear operators},  Wiley.
}

Honda, T. (2010), `Nonparametric regression for dependent data in the errors-in-variables problem`, Global COE Hi-Stat Discussion Paper Series, Institute of Economic Research, Hitotsubashi University.

\textcolor{blue}{
Hong, S. Y. and   Linton, O. B. (2016),   `Asymptotic properties of a Nadaraya-Watson type estimator for regression functions of infinite order', \emph{SSRN Electronic Journal}.
}

Jennen-Steinmetz, C. and Gasser, T. (1989), `A unifying approach for nonparametric regression estimation', \emph{J. Americ. Stat. Assoc.} 83, 1084--1089.

\textcolor{blue}{
 Jiang, J. and   Mack, Y.P. (2001), `Robust local polynomial regression for dependent data', \emph{Statistica Sinica} 11, 705--722.
}

Jones, M.C., Davies, S.J., and Park, B.U. (1994), `Versions of kernel-type regression estimators',
\emph{J. Americ. Stat. Assoc.} 89, 825--832.

\textcolor{blue}{
 Karlsen, H.A., Myklebust, T., and   Tjostheim, D. (2007), `Nonparametric estimation in a nonlinear
cointegration type model', \emph{Ann. Statist.} 35, 252--299.
}

Kulik, R. and Lorek, P. (2011), `Some results on random design regression with long memory errors and predictors', \emph{J. Statist. Plann. Infer.} 141, 508--523.

Kulik, R. and Wichelhaus C. (2011), `Nonparametric conditional variance and error density estimation in regrssion models with dependent errors and predictors', \emph{Electr. J. Statist.} 5, 856--898.

Laib, N. and Louani, D. (2010), `Nonparametric kernel regression estimation for stationary ergodic data: Asymptotic properties',
{\it J. Multivar. Anal.} 101, 2266--2281.

\textcolor{blue}{
Liang, H.-Y. and Jing, B.-Y. (2005),  `Asymptotic properties for estimates of nonparametric regression models based on negatively associated sequences', \emph{J. Multivariate Anal.} 95, 227--245.
}

Li, Y. and Hsing, T. (2010), `Uniform convergence rates for nonparametric
regression and principal component analysis in
functional/longitudinal data', \emph{Ann. Statist.} 38, 3321--3351.

\textcolor{blue}{
Li, X.,  Yang, W.,  and  Hu, S.  (2016), `Uniform convergence of estimator for nonparametric regression with dependent data', \emph{J. Inequal. Appl.},   142.
}

\textcolor{blue}{
 Lin, Z. and  Wang, J.-L. (2022), `Mean and covariance estimation for functional snippets', \emph{J. Amer. Statist. Assoc.} 117, 348--360.
}

\textcolor{blue}{
Linton,  O. and Wang, Q. (2016), `Nonparametric transformation regression with nonstationary data', \emph{Econometric Theory} 32, 1--29.
}

 Linke, Yu.Yu. and Borisov, I.S. (2017), `Constructing initial estimators in one-step estimation procedures of nonlinear regression',  \emph{Stat. Probab. Lett.}  120, 87--94.

 Linke, Yu.Yu. and Borisov, I.S. (2018), `Constructing explicit estimators in nonlinear regression problems',  \emph{Theory Probab. Appl.}  63, 22--44.

\textcolor{blue}{
 Linke, Yu.Yu. (2019), `Asymptotic properties of one-step  M-estimators', \emph{Commun. Stat. Theory Methods} 48, 4096--4118.
}

\textcolor{blue}{
{Linke, Yu.Yu. } (2023), `Towards insensitivity of Nadaraya--Watson estimators to design correlation', \emph{Theory Probab. Appl.} 68 (to appear).
}

\textcolor{blue}{
Linke, Yu.Yu. and Borisov, I.S. (2022), `Insensitivity of Nadaraya--Watson estimators to design correlation', \emph{Commun. Stat. Theory Methods} 51,  6909--6918.
}

\textcolor{blue}{
Linton, O. B. and   Jacho-Chavez, D. T. (2010), `On internally corrected and symmetrized kernel estimators for nonparametric regression', \emph{TEST} 19,  166--186.
}

\textcolor{blue}{
Loader, C. (1999),  \emph{Local regression and likelihood},  Springer.
}

Mack, Y.P. and M\"uller, H.-G. (1988), `Convolution type estimators for nonparametric regression', \emph{Stat. Prob. Lett.}  7, 229--239.

Masry, E. (2005), `Nonparametric regression estimation for dependent functional data', {\it Stoch. Proc. Their Appl.} 115, 155--177.

 M\"{u}ller, H.-G. (1988), \emph{Nonparametric Regression Analysis of Longitudinal Data}, New York: Springer.

\textcolor{blue}{
M{\"u}ller, H. G. and  Prewitt, K. A. (1993), `Multiparameter bandwidth processes and adaptive surface smoothing', \emph{J. Multivariate Anal.} 47, 1--21.
}

 Priestley, M.B. and Chao, M.T. (1972), `Non-Parametric Function Fitting',
 \emph{J. Royal Statist. Soc.}, Series B, 34, 385--392.

Rio, E. (2009), `Moment Inequalities for Sums of Dependent Random Variables under Projective Conditions', {\it J. Theor. Probab.}  22: 146--163.

Roussas, G.G. (1990), `Nonparametric regression estimation under mixing conditions', {\it Stach. Proc. Appl.}, 36, 107-116.

Roussas, G.G. (1991), `Kernel estimates under association: Strong uniform consistency', {\it Stat. Probab. Lett.} 12, 393-403.

\textcolor{blue}{
Roussas, G. G., Tran, L. T.,  and Ioannides, D. A. (1992), `Fixed design regression for time series: asymptotic normality', \emph{J. Multivariate Anal.}  40, 262--291.
}

\textcolor{blue}{
Shen, J. and   Xie,  Y. (2013),  `Strong consistency of the internal estimator of nonparametric regression with dependent data', \emph{Stat.  Probab. Lett.} 83, 1915--1925.
}

\textcolor{blue}{
Song, Q., Liu, R., Shao, Q., and  Yang, L. (2014), `A simultaneous confidence band for dense longitudinal regression', \emph{Commun. Stat. Theory Methods} 43, 5195--5210.
}

\textcolor{blue}{
Tang, X., Xi, M.,   Wu, Y.,  and Wang, X. (2018), `Asymptotic normality of a wavelet estimator for asymptotically negatively associated errors', \emph{Stat.  Probab. Lett.} 140, 191--201.
}

Wand, M.P. and Jones, M.C. (1995), \emph{Kernel Smoothing}, London: Chapman and Hall.

\textcolor{blue}{
Wang, J.-L., Chiou, J.-M., and M{\"u}ller, H.-G. (2016), `Review of functional data analysis', \emph{Annu. Rev. Statist.} 3, 257--295.
}

\textcolor{blue}{
Wang, Q. and  Chan, N. (2014), `Uniform convergence rates for a class of martingales with application in non-linear
cointegrating regression', \emph{Bernoulli} 20, 207--230.
}

\textcolor{blue}{
Wang, Q.Y. and   Phillips,  P.C.B. (2009a),  `Asymptotic theory for local time density estimation and nonparametric cointegrating regression', \emph{Econometric Theory} 25, 710--738.
}

\textcolor{blue}{
Wang, Q. and   Phillips, P.C.B. (2009b), `Structural nonparametric cointegrating regression', \emph{Econometrica} 77, 1901--1948.
}

Wu, J.S. and Chu, C.K. (1994), `Nonparametric estimation of a regression function with dependent observations',
{\it Stoch. Proc. Their Appl.}, 50, 149-160.

\textcolor{blue}{
Wu, Y., Wang, X., and  Balakrishnan, N. (2020), `On the consistency of the P--C estimator in a nonparametric regression model', \emph{Stat. Papers} 61, 899--915.
}

\textcolor{blue}{
Yang, X. and Yang, S. (2016), `Strong consistency of non parametric kernel regression estimator for strong mixing samples', \emph{Commun. Stat. Theory Methods}  46, 10537--10548.
}

\textcolor{blue}{
 Yao, F. (2007),  `Asymptotic distributions of nonparametric regression estimators for longitudinal or functional data', \emph{J. Multivariate Anal.}   98, 40--56.
}

\textcolor{blue}{
Yao,  F.,  M{\"u}ller, H.-G., and Wang, J.-L. (2005), `Functional data analysis for sparse longitudinal data', \emph{J. Amer. Statist. Assoc.} 100,  577--590.
}

\textcolor{blue}{
Young, D.S. (2017), \emph{Handbook of regression methods}, Chapman and Hall.
}

\textcolor{blue}{
Zhang, X. and  Wang,  J.-L. (2016),  `From sparse to dense functional data and beyond', \emph{Ann. Statist.} 44, 2281--2321.
}

\textcolor{blue}{
Zhang, J.-T. and  Chen, J. (2007), `Statistical inferences for functional data', \emph{Ann. Statist.} 35, 1052--1079.
}

\textcolor{blue}{
Zhang, X. and   Wang, J.-L. (2018), `Optimal weighting schemes for longitudinal and functional data', \emph{Stat.  Prob. Lett.} 138, 165--170.
}

\textcolor{blue}{
Zhang, S., Miao, Y., Xu, X., and   Gao, Q. (2018),  `Limit behaviors of the estimator of nonparametric regression model based on martingale difference errors', \emph{J. Korean Stat. Soc.} 47, 537--547.
}

\textcolor{blue}{
Zhang, S., Hou, T.,  and Qu, C. (2019), `Complete consistency for the estimator of nonparametric regression model based on martingale difference errors', \emph{Commun. Stat. Theory Methods} 50, 358--370.
}

\textcolor{blue}{
Zhou, X. and Zhu, F. (2020), `Asymptotics for L1-wavelet method for nonparametric regression', \emph{J. Inequal. Appl.} 216.
}

\end{document}